\documentclass[11pt,reqno]{amsart}

\usepackage{graphicx,color}
\usepackage{amsmath}
\usepackage{amssymb}
\usepackage{amsfonts}
\usepackage{amsthm}
\usepackage{epstopdf}
\usepackage[color]{showkeys}
\newcommand{\R}{\mathbb R}
\newcommand{\Z}{\mathbb Z}
\newcommand{\T}{\mathbb T}

\newtheorem{conjecture}{Conjecture} 
\theoremstyle{remark}
\newtheorem{remark}{Remark}[section]
\theoremstyle{definition}
\newtheorem{definition}{Definition}[section]

\begin{document}
\title[Line solitons of the 2D Zakharov-Kuznetsov equation]
{Numerical study of the transverse stability of line solitons of the Zakharov-Kuznetsov equations}

\author[C. Klein]{Christian Klein}
\address{Institut de Math\'ematiques de Bourgogne,  UMR 5584;\\
Institut Universitaire de France \\
Universit\'e de Bourgogne-Franche-Comt\'e, 9 avenue Alain Savary, 21078 Dijon
                Cedex, France}
\email{Christian.Klein@u-bourgogne.fr}

\author[J.C. Saut]{Jean-Claude Saut}
\address{Laboratoire de Mathématiques\& Universit\'e de Paris - Saclay,  
91405 Orsay, France.}
\email{jean-claude.saut@universite-paris-saclay.fr}
    
\author[N. Stoilov]{Nikola Stoilov}
\address{
Laboratoire Jaques-Louis Lions, UMR 7598, Sorbonne Universit\'e, 4 Place Jussieu
75005 Paris, France.
} 
\email{Nikola.Stoilov@ljll.math.upmc.fr}
\begin{abstract}
We present a detailed numerical study of the stability 
under periodic perturbations of 
line solitons of two-dimensional, generalized Zakharov-Kuznetsov equations 
with various power nonlinearities. In the 
$L^{2}$-subcritical case, in accordance with a theorem due to 
Yamazaki we find a critical speed, below which the line soliton 
is stable. For higher velocities, the numerical results indicate an 
instability against the formation of lumps, solitons localized in both 
spatial directions.  In the $L^2$-critical and supercritical cases 
but subcritical for the 1D generalized Korteweg-de 
Vries equation), 
the line solitons are shown to be numerically stable for small 
velocities, and strongly unstable for large velocities, 
with a blow-up observed in finite time. 
\end{abstract}

\thanks{This work was supported by the ANR-FWF project ANuI - 
ANR-17-CE40-0035. CK and NS thank for support by  the isite BFC project NAANoD, the ANR-17-EURE-0002 EIPHI and by the European Union Horizon 2020 research and innovation program under the Marie Sklodowska-Curie RISE 2017 grant agreement no. 778010 IPaDEGAN.}

\maketitle
\section{Introduction}

This paper is concerned with the stability of line solitons to the two-dimensional (2D) generalized Zakharov-Kuznetsov 
(ZK) equation
\begin{equation}
    u_t + (u_{xx}+u_{yy} + u^p)_x  = 0,
    \label{ZK}
\end{equation}
or perturbations  periodic in $y$.
This equation is an extension of the generalized 
Korteweg-de Vries (KdV) to two spatial dimensions. 
It is not necessary for $p$ to be an integer in this equation, however in this paper we 
will only consider integer nonlinearities, concretely the cases 
$p=2,3,4$. 
The  ZK 
equation in 2D with a quadratic nonlinearity ($p=2$) was originally 
proposed by Zakharov and Kuznetsov   \cite{ZK}, see also \cite{KRZ}, and  rigorously 
justified in  \cite{LLS} (see also \cite{Pu}) from the Euler-Poisson 
system  for  uniformly magnetized plasma. It was also rigorously 
derived by Han-Kwan (\cite{Han-Kwan}) from the Vlasov-Poisson system 
in the presence of an external magnetic field. Note that in those two physical contexts, the presence of an applied magnetic field explains the lack of spatial symmetry of the ZK equation.

The generalized ZK equation with $p=3$ (`modified ZK equation')
 appears as  an asymptotic model in  the context of weakly 
nonlinear ion-acoustic waves in a plasma of cold ions and hot 
isothermal electrons with a uniform magnetic field 
\cite{MP1999}. It can be also applied  as the 
amplitude equation for two-dimensional long waves on the free surface 
of a thin film flowing down a vertical plane with moderate values of 
the fluid surface tension and large viscosity, see \cite{MM1989}. 
On the other hand, the case $p=4$ does not seem to  appear as a physically relevant model but (as the generalized KdV equation) can 
be used as a mathematical toy model to investigate the competition between nonlinearity and dispersion.

The generalized KdV equations have solitary wave solutions of 
the form $u=Q_{c}(z)$ with $z=x-x_{0}-ct$, $x_{0},c=const$, 
$x_{0}\in\mathbb{R}$, 
$c>0$, and with
\begin{equation}
    Q_{c}(z) = \left( 
    \frac{(p+1)c}{2}\,\mbox{sech}^{2}\frac{\sqrt{c}(p-1)}{2}z\right)^{1/(p-1)}
    \label{soliton}.
\end{equation}

It is well-known that these solitary waves are 
orbitally  and asymptotically stable in the context of the 
generalized KdV equation in the $L^2-$ subcritical case  $p< 4$.

These solitary waves can form $y$-independent solutions of the ZK 
equations by a trivial extension in the $y$-direction and are called 
\emph{line solitons}. This notion can obviously be also 
applied to settings periodic in $y$ as exclusively considered here. A natural question is that of the {\it 
transverse stability } of the line  solitons which will be studied in 
two different settings, depending of the nature of the perturbation, 
either for localized (in $(x,y)$ perturbations of the line soliton, 
or for $y$ periodic perturbations). We refer to Def. 
\ref{definition} for a precise definition of the used concept of stability.

In this paper we study numerically the transverse stability of these 
line solitons for power nonlinearities $p=2,3,4$. This extends the 
numerical study of \cite{KRS} (see \cite{KRS2} for 3D simulations) 
for localized initial data to solutions being 
localised in $x$ and periodic in $y$. Note 
that the 2D ZK equations also have solitary wave solutions 
$\mathcal{Q}_{c}$ to be defined in (\ref{E:Q}), (\ref{Q}) and 
(\ref{Qscal}) which are 
exponentially localized in all spatial directions called 
\emph{lumps}, see \cite{dB} for  stability of such 
solutions, \cite{CMPS} for their asymptotic stability  and  
\cite{KRS} for figures.   

We briefly recall some relevant theoretical results.
The first transverse instability result of the  KdV line soliton with respect to (localized) two-dimensional perturbations was given by Rousset and Tzvetkov  in \cite{RT2} in the spirit of similar results (\cite{RT})  for the nonlinear stability of the KdV soliton 
as a line soliton of the Kadomtsev-Petviashvili (KP) equation.
Bridges \cite{Bridges}  showed the  
instability of the line solitary waves $Q_c( x-ct)$ of the 
Zakharov-Kuznetsov equation on $\R\times \T_L$ with sufficiently 
large speed $c$ (where $\T_L$ is the torus with $2\pi 
L$ period, that is $\T_L=\R/(2\pi L \Z)$), see also \cite{Kuz}.

In the case $p=2$, Yamazaki \cite{YAM17} applying the approach by Rousset 
and Tzvetkov \cite{RT} proved that the line soliton of the 
Zakharov-Kuznetsov equation on $\R\times \T_L$ 
 is orbitally,  and 
moreover asymptotically stable for $0<c\leq c^{*}:=\frac{4}{5L^2}$ 
and is orbitally unstable  for $c>\frac{4}{5L^2}.$ His definition of 
orbital stability reads as follows:
\begin{definition}\label{definition}
	We say that a line solitary wave $Q_{c}(x - ct, 
	y)$ is orbitally transversally stable in 
	$\mathbb{H}^{1}(\mathbb{R} \times \mathbb{T}_{L})$ if for any 
	$\epsilon>0$ there exists $\delta>0$ such that for all initial 
	data $u_{0}\in\mathbb{H}^{1}(\mathbb{R} \times \mathbb{T}_{L})$ 
	with
	$|u_{0}-Q_{c}|_{H^{1}} <\delta$, the solution $u(t )$ of (\ref{ZK}) 
	with $u(0) = u_{0} $ exists globally in positive time and satisfies
$$\mbox{sup}_{t>0} \mbox{inf}_{(x_{0},y_{0})\in \mathbb{R} \times \mathbb{T}_{L}} 
|u(t,\cdot,\cdot)-Q_{c}(\cdot-x_{0},\cdot-y_{0})|_{H^{1}}<\epsilon.$$
Otherwise, we say the solitary wave $Q_{c}(x - ct, y)$ is orbitally 
transversally unstable in $\mathbb{H}^{1}(\mathbb{R} \times \mathbb{T}_{L})$.
\end{definition}
We apply this definition of stability throughout the paper, so that we are only concerned with transverse stability (instability), the perturbation being always $y$-periodic and not localized in $(x ,y)$. 
Stability with respect to perturbations with a $y$-dependence will be 
referred to as \emph{transverse stability.}

This result was completed in \cite{YAM20} where Yamazaki  was able to construct center stable manifolds around unstable line solitary waves to the Zakharov-Kuznetsov equation on $\R\times \T_L.$
Pelinovsky \cite{Pel} proved the asymptotic stability of the transversely modulated solitary waves of the Zakharov-Kuznetsov equation on $\R\times \T_L$ in exponentially weighted spaces.

Note that the type of the instability for $c>c^{*}$ is unknown, and we state 
 the following  conjecture based on our numerical experiments\\
\textbf{Main conjecture I}:\\
\emph{Consider equation (\ref{ZK}) for $(x,y)\in \mathbb{R}\times 
\mathbb{T}_{L}$. The line solitons (\ref{soliton}) for $p=2$ are for 
$c>c^{*}(L)$ 
unstable against the formation of lumps.}

This conjecture was first suggested by the pioneering numerical simulations in \cite{ITK}, and is similar to what is found numerically in the context 
of the KP I equation, see \cite{KS12}. 

In the cases $p=3,4$ for the ZK equation, the line solitons are strongly unstable 
(i.e., leading to a blow-up in finite time of the $L^{\infty}$ norm 
of the solution), such as the 
corresponding lumps studied in \cite{KRS}, i.e., the line soliton is 
unstable against lumps, but the latter then blow up in finite time. 
We get \\
\textbf{Main conjecture II}:\\
\emph{Consider equation (\ref{ZK}) for $(x,y)\in \mathbb{R}\times 
\mathbb{T}_{L}$. For $p=3,4$ the line solitons are strongly unstable. Perturbations 
with $c$ smaller than some critical speed $c^{*}(p,L)$  are  dispersed, 
perturbations with $c>c^{*}(p,L)$ lead to a blow-up in finite time. The 
blow-up mechanism  (i.e., blow-up rate and profile) is as conjectured in \cite{KRS} for a situation in 
$\mathbb{R}^{2}$, see conjectures 
\ref{C:2}, \ref{C:3}.
}

The paper is organized as follows. In section 2 we collect some basic 
facts on the ZK equation and the used numerical tools. In section 3 
we consider the subcritical case $p=2$, both in the stable and in the 
unstable regime. In section 4 we present 
numerical results for the critical case and show that a blow-up can 
be observed in some cases. A similar study for the supercritical case 
$p=4$ is presented in section 5. We add some concluding remarks in 
section 6. 

\section{Basic facts}
In this section we collect some basic facts on the ZK and generalized ZK equations in 2D 
and on the used numerical approaches. 

\subsection{Analytic facts}

First, one notes that if $u$ solves the generalized ZK equation with 
initial data $u_0$,  then $u_\lambda(x,y,t)=\lambda^{2/(p-1)}u(\lambda 
x,\lambda y, \lambda ^3 t)$ is also a solution with initial data 
$u_{0,\lambda}(x,y)= \lambda^{2/(p-1)}u_{0}(\lambda x,\lambda y)$ for any $\lambda >0.$
This implies that the homogeneous Sobolev space $\dot{H}^s(\R^2)$ is invariant by this scaling  when $s=s_c=1-\frac{2}{p-1},$ in particular  $p=3$ is the $L^2-$ critical exponent.

We will not describe in details the many papers devoted to the 
well-posedness for the ZK equation  in the whole space $\R^2,$ starting with the pioneering work of Faminskii \cite{F95} who proved global well-posedness in $H^1(\R^2)$. The best local well-posedness result  is established by Kinoshita in \cite{K2021} in $H^s(\R^2), s>-1/4$, and this implies the global well-posedness in $L^2(\R^2)$. We refer to this last paper for an extensive bibliography.

As for the generalized  ZK equation when $p=3,4,$ Ribaud and Vento 
\cite{RV2013} proved local well-posedness in $H^s(\R^2)$, $s>1/4$ $(p=3)$, $s>5/12$ $(p=4)$. A finite time blow-up is  expected in these 
cases, possible consequence of the instability of the lump solitary 
wave proven in \cite{FHR3}. Actually, it was proven in \cite{FHRY} 
that a blow-up may occur in finite or infinite time in the cubic case 
($p=3).$ Such a blow-up result is expected but still  unproven in the 
$L^2$ supercritical case, $p=4.$ \footnote{Note that  finite time 
blow-up is still unproven for the general $L^2$ supercritical generalized KdV 
equation except in the perturbative case $p=5^{+}$ in \cite{Lan}.}

On the other hand, as recalled in the Introduction, the study of the transverse stability of line solitons can be performed in two 
settings depending on the nature of perturbations: either fully 
localized two-dimensional perturbations or periodic in $y$ perturbations.
A first step is to prove the  well-posedness of the Cauchy problem in those two settings.

We first focus on the usual ZK equation, $p=2.$
In this case, one has to solve the problem

\begin{equation}\label{IVP}
     \left\{
\begin{array}{lll}
{\displaystyle u_t+\partial_x \Delta u+u\partial_xu+\partial_x(\phi 
u)  =  0,  }  \qquad (x,y) \in \mathbb{R}^2, \,\,\,\, t>0, \\
{\displaystyle  u(x,y,0)=u_0(x,y)},
\end{array}
\right.
\end{equation}

where $\phi$ can be the KdV solitary wave or any KdV N-soliton. It was proven in \cite{LPS} that the Cauchy problem \eqref{IVP} is globally well-posed in $H^1(\R^2).$

The second situation necessitates to solve the Cauchy problem for the 
ZK equation in the spatial domain $\R\times \T.$ The local well-posedness in $H^s(\R\times \T),\; s>3/2$ is proven in \cite{LPS} and the global well-posedness in $H^1(\R\times \T)$ is proven in \cite{MoPi}.

We are not aware of similar results for the generalized ZK equation 
($p=3,4)$.

In addition to the line solitons which are just $y$-independent 
solutions of the KdV equations (\ref{soliton}), the 2D ZK equation 
has a family of traveling wave solutions  localized in both spatial 
directions called \emph{lumps}
\begin{equation}\label{E:Q}
u(x,y,t) = \mathcal{Q}(x-ct,y)
\end{equation}
satisfying
\begin{equation}\label{Q}
-c\mathcal{Q}+\mathcal{Q}_{xx}+\mathcal{Q}_{yy}+\mathcal{Q}^{p}=0;
\end{equation}
the solitary waves \( \mathcal{Q}_{c}(x,y) \) are related to  \(\mathcal{Q}_{1}(x,y)=:\mathcal{Q}(x,y) \) 
for \( c > 0 \) via
\begin{equation} \label{Qscal}
\mathcal{Q}_{c}(x,y)=c^{\frac1{p-1}} 
\,\mathcal{Q}(\sqrt{c}\,x,\sqrt{c}\,y).
\end{equation}

As noticed in \cite{dB} the existence of ground state solutions, that is positive, radially symmetric solutions, results from the work of Berestycki and Lions \cite{BL} see also \cite{SS1999}. Moreover  $\mathcal{Q}$ is smooth and decays exponentially at infinity. 
We refer for instance to \cite{ITK, KRS} for numerical simulations. 
The ground state is unique (up to standard symmetries) by a classical 
result of Kwong \cite{Kw}. A. de Bouard \cite{dB} proved its orbital 
stability, the asymptotic stability in the subcritical case was 
proven in \cite{CMPS}.

The above scaling invariance with respect to $\lambda$ for the equation \eqref{ZK},
can be used
in the context of blow-up in the form of a \emph{dynamical rescaling}
\begin{equation}\label{resc}
\begin{array}{c}
X = \frac{x-x_{m}(t)}{\lambda(t)}, \quad Y = \frac{y-y_{m}(t)}{\lambda(t)}, \quad 
    T=\int_{0}^{t}\frac{dt'}{\lambda^{3}(t')},\\
    ~\\
U(X,Y,T) = \lambda(t)^{\frac2{p-1}}(t) \, u(x,y,t).
\end{array}
\end{equation}
The dynamically rescaled ZK equation reads
\begin{equation}\label{ZKresc}
U_T - a\bigg(\frac2{p-1} U+X U_{X}+Y U_{Y}\bigg) - v_{X}U_{X}-v_{Y}U_{Y}+ \bigg(U_{XX}+U_{YY} + U^p\bigg)_X  = 0,
\end{equation}
where 
\begin{equation}\label{a}
a \equiv a(T) = \frac{d\ln \lambda}{dT},\quad v_{X}=\frac{x_{m,T}}{\lambda},\quad v_{Y}=\frac{y_{m,T}}{\lambda}.
\end{equation}
A potential blow-up is expected for $T\to\infty$, where \( 
U_{T} \) is assumed to vanish. Thus, the equation 
\eqref{ZKresc} in the limit becomes 
\begin{equation}\label{ZKresinfty}
-\overset{\infty}{a}\bigg(\frac{2}{p-1}\overset{\infty}{U}+ X\overset{\infty}{U}_X+Y\overset{\infty}{U}_Y \bigg) - 
 v_{\underset{\infty}{X}}\overset{\infty}{U}_{X}-v_{\underset{\infty}{Y}}\overset{\infty}{U}_{Y}+ \bigg(\overset{\infty}{U}_{XX}+\overset{\infty}{U}_{YY} + \overset{\infty}{U}~^p \bigg)_X = 0,
\end{equation}
where the sub/superscript \( \infty \) denotes that the quantity is taken in the limit as $T\to\infty$ and $ \overset{\infty}{U} $ stands for a blow-up profile. 

As discussed in \cite{KRS} two possible stable blow-up mechanisms are expected: either an algebraic dependence of \( \lambda \) on \( T \), or an exponential 
one. In the former case the quantity \( \overset{\infty}{a} \) in 
\eqref{a} will vanish, and equation \eqref{ZKresinfty} will be 
identical to the equation for the lump if 
$v_{\underset{\infty}{Y}}=0$;
this mechanism is expected in the $L^{2}$-critical case. If 
$\lambda\propto 1/T$ as in the $L^{2}$-critical generalized KdV case, one 
has
\begin{equation}\label{Lcritical}
\lambda\propto \sqrt{t^{*}-t}.
\end{equation}
In the supercritical case, one expects an exponential decay of $\lambda$ 
with $T$, that is, $\lambda\propto \exp(-\gamma T)$ with $\gamma>0$, 
\begin{equation}\label{Lsuper}
\lambda\propto (t^{*}-t)^{1/3}.
\end{equation}

In \cite{KRS}, numerical results led to the following conjectures
\begin{conjecture}[$L^2$-critical case] \label{C:2} 
Consider the critical 2D ZK equation \eqref{ZK} with $p=3$.
\begin{enumerate}
\item
If $u_0 \in \mathcal{S}(\mathbb{R}^{2})$ is such that 
$\|u_0\|_{2} < \|\mathcal{Q}\|_{2}$, then the solution $u(t)$ to \eqref{ZK} is dispersed. 
\item    
If $u_0 \in \mathcal{S}(\mathbb{R}^{2})$ is sufficiently localized and such that 
$\|u_0\|_{2} > \|\mathcal{Q}\|_{2}$, then the solution blows up in finite time $t=t^*$ and
such that as $t \to t^*$
\begin{equation}\label{selfs}
u(x,y,t)- \frac{1}{\lambda(t)}\,\mathcal{Q}\left(\frac{x-x_{m}(t)}{\lambda(t)}, \frac{y-y_{m}(t)}{\lambda(t)} \right) \to \tilde{u}\in L^{2},
\end{equation}
with 
\begin{equation}\label{ux2}
\|u_{x}(t)\|_{2} \sim \frac{1}{\lambda(t)}, ~~ \lambda(t)\sim \sqrt{t^*-t}, \quad \mbox{and} \quad 
x_{m}(t) \sim \frac{1}{t^*-t}, ~~ y_m(t) \to y^* <\infty.
\end{equation}
\end{enumerate}
\end{conjecture}

\begin{conjecture}[$L^2$-supercritical case] \label{C:3}
Consider the supercritical 2D ZK equation, in particular, when $p=4$ in \eqref{ZK}.
Let $u_0 \in \mathcal{S}(\mathbb R^2)$ be of sufficiently large mass and energy 
 and of some localization. 
Then the ZK solution $u(t)$ blows up in finite time $t^{*}$ and finite location $(x^{*},y^{*})$, i.e., the blow-up core resembles a self-similar structure with  
\begin{equation}\label{selfss}
u(x,y,t)- \frac{1}{\lambda^{\frac2{p-1}}(t)}\, P\left(\frac{x-x_{m}(t)}{\lambda(t)}, \frac{y-y_{m}(t)}{\lambda(t)}\right)\to \tilde{u}\in L^{2},
\end{equation}
where $P(x,y)$ is a localized solution to \eqref{ZKresinfty} (which 
is conjectured to exist),  
$$
x_m(t) \to x^*, \quad y_m(t) \to y^{*},
$$
and 
\begin{equation}\label{ux2s}
\|u_{x}(t)\|_{2}\sim \frac{1}{\lambda^{\frac2{p-1}}(t)} \quad \mbox{with} \quad  \lambda(t)\sim (t^{*}-t)^{1/3} \quad \mbox{as} \quad t \to t^*.
\end{equation}
\end{conjecture}

\subsection{Numerical approaches}
The numerical approach for ZK is as in \cite{KRS} to which the reader 
is referred for details. Though we study an analytical situation 
on $\mathbb{R}\times \mathbb{T}_{L}$, we approximate this numerically by 
working on $\mathbb{T}_{L_{x}}\times\mathbb{T}_{L_{y}}$. Thus both the $x$ and $y$ 
dependence is approximated via a discrete Fourier transform 
implemented with a \emph{fast Fourier transform}. This can be seen as 
a truncated Fourier series. Since it is well known that the Fourier 
coefficients decrease exponentially with the index for analytic 
functions, the numerical error in truncating decreases in the same 
way with the numerical resolution (i.e., the number $N_{x}$, $N_{y}$ 
of Fourier modes in $x$ respectively $y$). Therefore the Fourier 
coefficients can be used to indicate the numerical resolution 
in the computational domain. 

Since we want to approximate a situation on $\mathbb{R}\times 
\mathbb{T}_{L}$ by a model in $\mathbb{T}^{2}$, we have to choose the 
period in $x$ large enough that radiation emitted towards $-\infty$ 
does not have a significant effect on the studied phenomena. The 
period $2\pi L_{x}$ in $x$ is always  
chosen  sufficiently large so that the line 
soliton decreases to the order of machine precision (roughly $10^{-15}$ in 
double precision) and that no significant amount of radiation re-enters the 
computational domain during the computation, i.e., that the amplitude 
of the radiation at the location of the line soliton is much smaller than the 
amplitude of the latter. 

The integration in time is carried out with an \emph{exponential time 
differencing scheme} by Cox and Matthews \cite{CM} since these were 
the most efficient for KdV equations in \cite{etna,KR}, see \cite{HO} for 
a review on such integrators. Since the mass and the energy, exactly 
conserved quantities of the ZK equation, are not automatically 
conserved by the time integration scheme, they can be used as 
discussed in \cite{etna,KR} to control the resolution in time for a 
given resolution in space (typically the relative conservation of 
these quantities overestimates the accuracy in an $L^{\infty}$ sense 
by 2-3 orders of magnitude). In this paper we will always use the 
relative conservation of the mass to indicate the accuracy in the 
time integration. Note that here the mass is considered on 
$\mathbb{T}^{2}$, not in $\mathbb{R}^{2}$ as usual which implies that 
the line soliton in this setting has finite mass. 

In this paper, we consider initial data representing perturbed line 
solitons. In particular we study 
 localized perturbations of the line soliton (\ref{soliton}),
\begin{equation}
u_0(x,y) = Q_{c}(x) + a e^{-b^2(x^2 +y^2)},  
\label{initial_gauss}
\end{equation} 
and a deformation of the soliton which is just localized in $x$, 
but slowly modulated in the $y$ direction of the form 
\begin{equation}
u_0(x,y) = Q_{c}(x) + a\cos^2(y/b+\delta)e^{-x^2},  
\label{initial_periodic}
\end{equation} 
where $a$, $b$, and $\delta$ are constants. 

Note that we do not solve the dynamically rescaled ZK equation 
(\ref{ZKresc}) even when a blow-up is observed. As in \cite{KRS} we 
always solve (\ref{ZK}) in such cases in order to avoid problems at 
the computational boundary (see the discussion in \cite{KP2}). The 
blow-up rate, i.e., $\lambda(\tau)$, is obtained by tracing the $L^{\infty}$ norm 
of the solution and some post-processing as in \cite{KRS}. 

\section{Sub-critical case}
In this section we study perturbations of the KdV soliton for the 2D 
ZK equation, i.e., the subcritical case. We consider both stable 
cases and an instability against lump formation.

\subsection{Stable regime}
As discussed above, we choose a large value for the period in $x$ in 
order to delimit the amount of radiation reentering the computational 
domain, here  $L_{x}= 250$ and $L_y=1$, 
which gives the critical speed $c^* = 4/5$. We discretise the domain 
by taking $N_{x}=2^{15}$ points in the $x$-direction and $N_{y}=2^7$ in the 
transversal $y$-direction. High-index Fourier coefficients, which are 
used to estimate the spatial resolution, stay below machine precision 
throughout the run, and the relative conservation of mass $\Delta = 
|1 - 
m/m_0|$ is below $10^{-12}$. 

First we consider localized perturbations of the line soliton of the 
form (\ref{initial_gauss}) with a 
speed  $c=0.75<c^*=4/5$; we choose $a = 0.1$  so that $a< 
0.1|u_0|_\infty$, and $b = 3$  in such a way that the 
perturbation is localized on the computational domain. 

\begin{figure}[!htb]
\includegraphics[width=0.51\hsize]{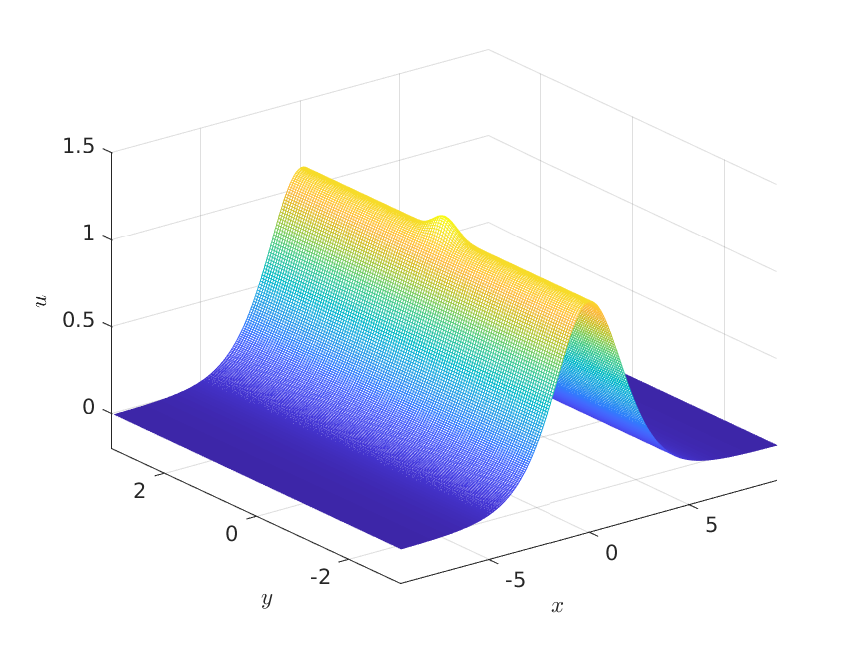}
\includegraphics[width=0.46\hsize]{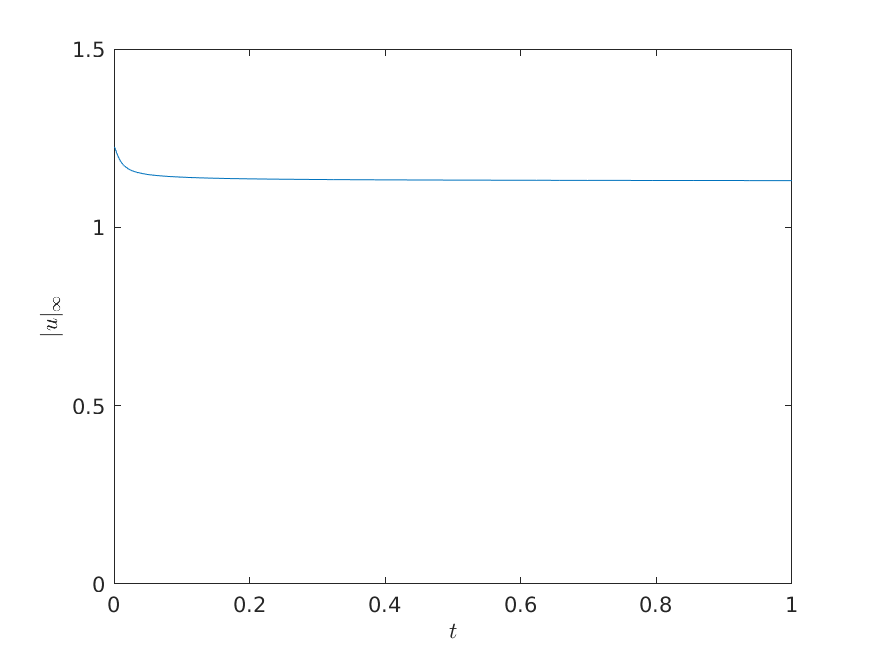} 
\caption{Initial data of the form 
(\ref{initial_gauss}) $u_0 = Q_c + a\exp^{-b^2r^2}$ with $a = 0.1$, 
$b = 3$  and sub-critical speed $c = 0.75$  for subcritical $(p =2)$ ZK in 
the stable regime on the left, 
and the time evolution of $|u|_\infty$ on the right.  
}
\label{fig:sub_stable_local_initial}
\end{figure}

The initial data can be seen on 
the left of Fig.~\ref{fig:sub_stable_local_initial}. For the time 
evolution we use $N_{t}=2000$ time steps for $t\in[0,1]$.  It can be seen on 
the right of Fig.~\ref{fig:sub_stable_local_initial} that   some radiation is emitted 
but that   the $L^{\infty}$ norm of the solution stabilises on some 
plateau.
   
\begin{figure}[!htb]
\includegraphics[width=0.49\hsize]{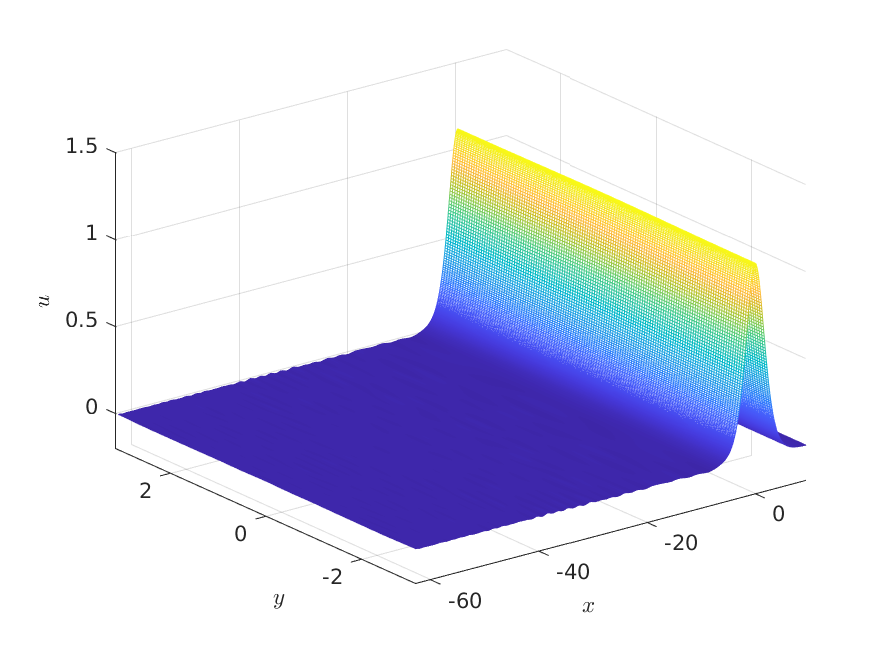} 
 \includegraphics[width=0.49\hsize]{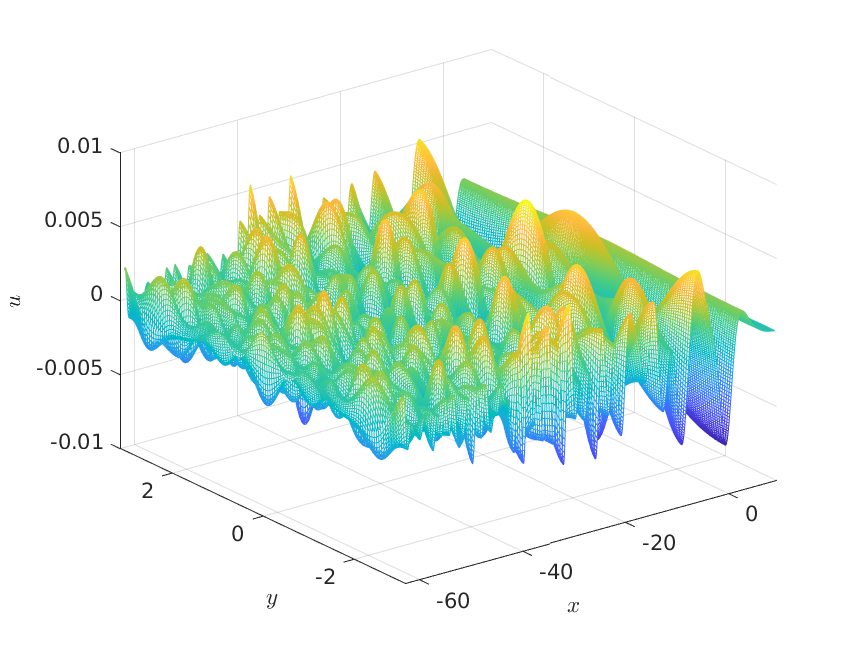}
\caption{Solution to the subcritical $(p =2)$ ZK for 
locally perturbed soliton initial data (\ref{initial_gauss}) $u_0 = Q_c + a\exp^{-b^2r^2}$ with $a = 0.1$, $b = 3$ 
 and sub-critical speed $c = 0.75$: on 
the left the solution for $t=1$  and on the right the difference between the final state and a 
fitted line soliton (\ref{soliton}) for the same time. The 
orientation of the plot is chosen to give clear view of the region behind the line soliton, where, in the stable regime, the interesting phenomena are located.} 
\label{fig:sub_stable_local_final}
\end{figure}

The behavior of the $L^{\infty}$ norm indicates that the final state 
is a line soliton with slightly different mass (in order to see 
numerically in finite time effects of the perturbation, we had to 
consider a perturbation with a mass of a few percent of the 
unperturbed soliton; thus though the soliton is clearly stable, the 
mass of the final state is slightly different from the one of the 
unperturbed soliton). On the left of Fig.~\ref{fig:sub_stable_local_final} we show 
the solution for $t=1$. To compare this final state 
to (\ref{soliton}), we simply determine the speed $c_{F}$ of the 
soliton (\ref{soliton}) having the same maximum as the final state. 
On the right of Fig.~\ref{fig:sub_stable_local_final} one can see that
the difference of the final state  and a line soliton with a fitted 
value $c_{F}$ of the speed  is within $0.5\%$ of the maximum of the 
initial data,  the same order of 
magnitude as the 
radiation. The fitted soliton has speed 
ratio $c_F/c = 0.9953$.

\begin{figure}[!htb]
\includegraphics[width=0.51\hsize]{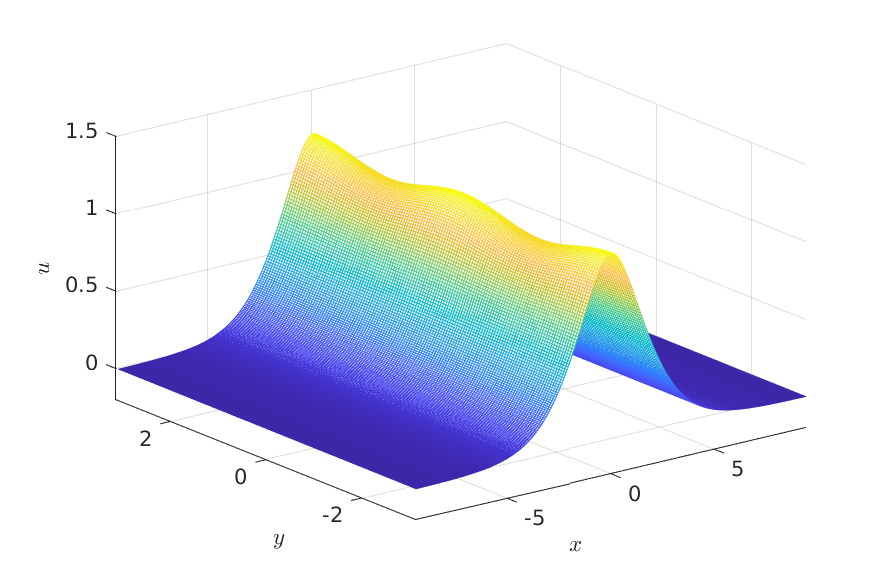}
\includegraphics[width=0.46\hsize]{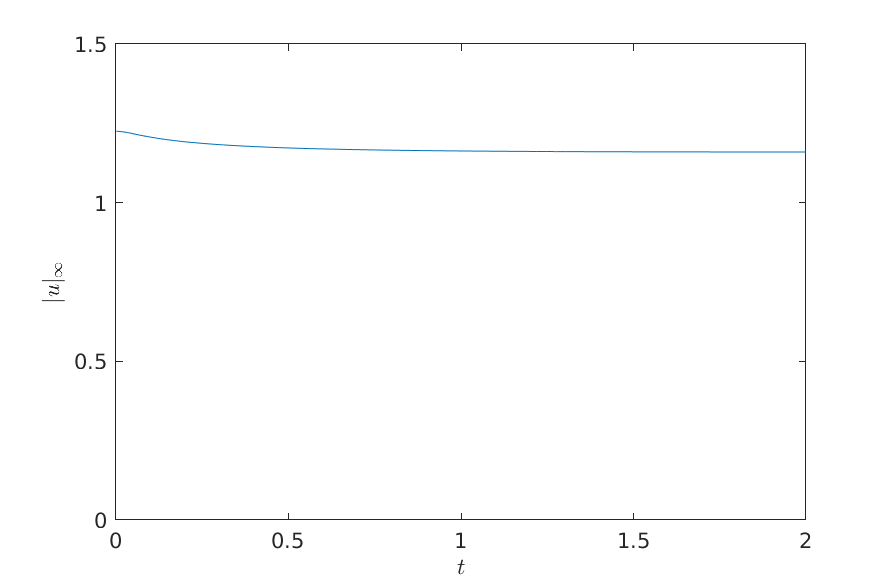} 
\caption{Initial data of the form 
(\ref{initial_periodic}) $u_0 = Q_c + a\cos^2(y/L_y)\exp^{-b^2x^2}$ with $a = 0.2$, $b = 1$  and sub-critical speed $c = 0.75$ on the left, and on the right time evolution of $|u|_\infty$ for subcritical $(p =2)$ ZK.}
\label{fig:sub_stable_periodic_initial}
\end{figure}

As a second example we consider non-localised, but periodic perturbations 
in $y$ of the form (\ref{initial_periodic}) with $c = 0.75$ and $a = 
0.1$, $b=1$, $\delta=0$ following the same reasoning as before. The initial data can be 
seen on the left of Fig.~\ref{fig:sub_stable_periodic_initial}. The time dependence of 
the $L^{\infty}$ norm on the right of the same figure once more 
indicates that the solution is stable, and that the final state is 
 a line soliton of slightly different mass.
             
\begin{figure}[!htb]
\includegraphics[width=0.49\hsize]{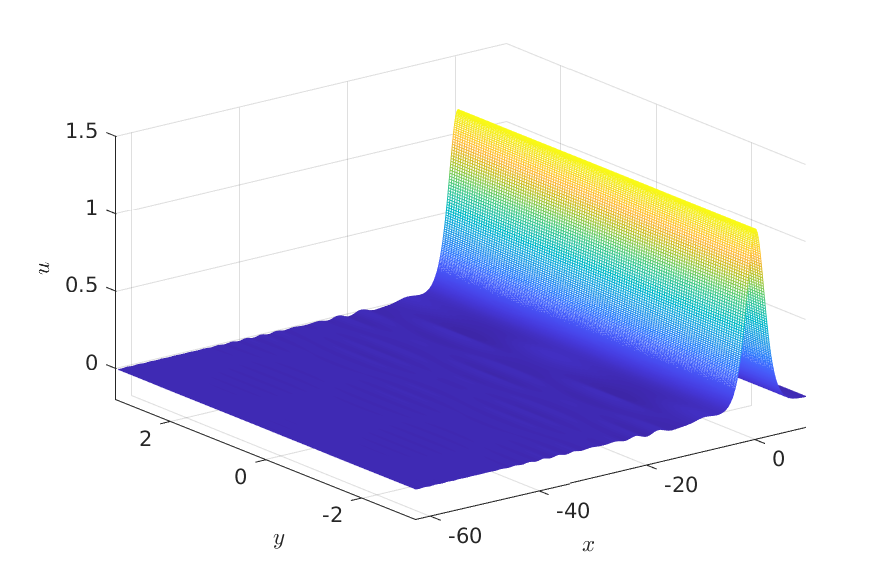} 
 \includegraphics[width=0.49\hsize]{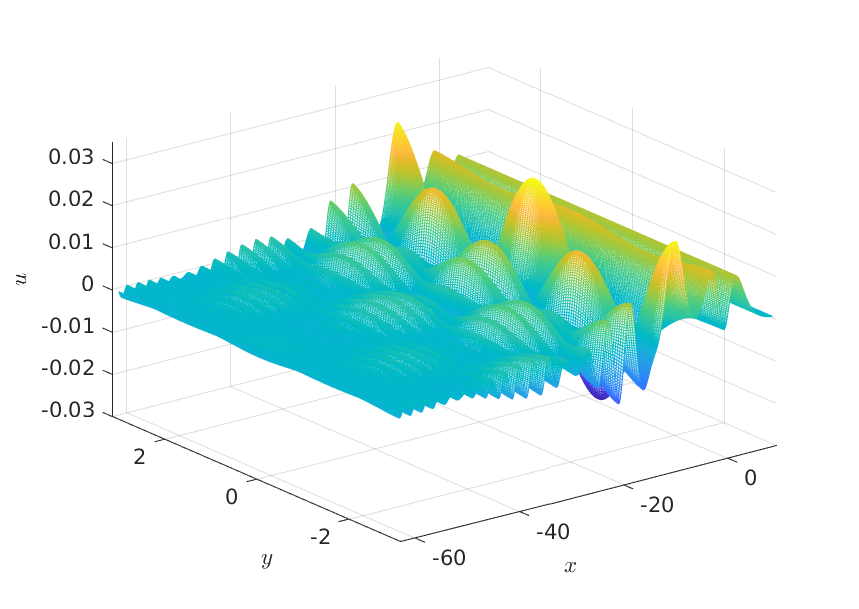}
\caption{Solution to the subcritical $(p =2)$ ZK for 
locally perturbed soliton initial data (\ref{initial_gauss}) $u_0 = Q_c + a\cos^2(y/Ly)\exp^{-b^2x^2}$ with $a = 0.2$, $b = 1$ 
 and sub-critical speed $c = 0.75$ at $t=1$ on 
the left, and on the right the difference between the final state and a 
fitted line soliton (\ref{soliton}) for the same time. }
\label{fig:sub_stable_periodic_final}
\end{figure}

The final state of the solution is shown on the left of 
Fig.~\ref{fig:sub_stable_periodic_final}. It is within $0.25\%$ of the fitted line 
soliton (\ref{soliton}) as can be seen on the right of the same 
figure. The fitted soliton has speed ratio $c_F/c = 0.9702$.

\subsection{Unstable regime}

\begin{figure}[!htb]
 \includegraphics[width=0.49\hsize]{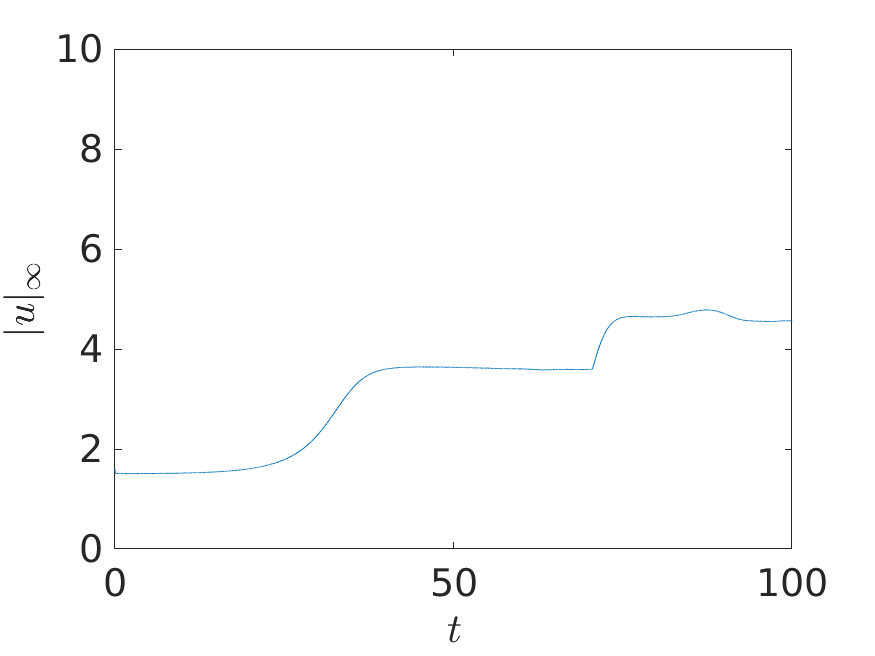} 
 \includegraphics[width=0.49\hsize]{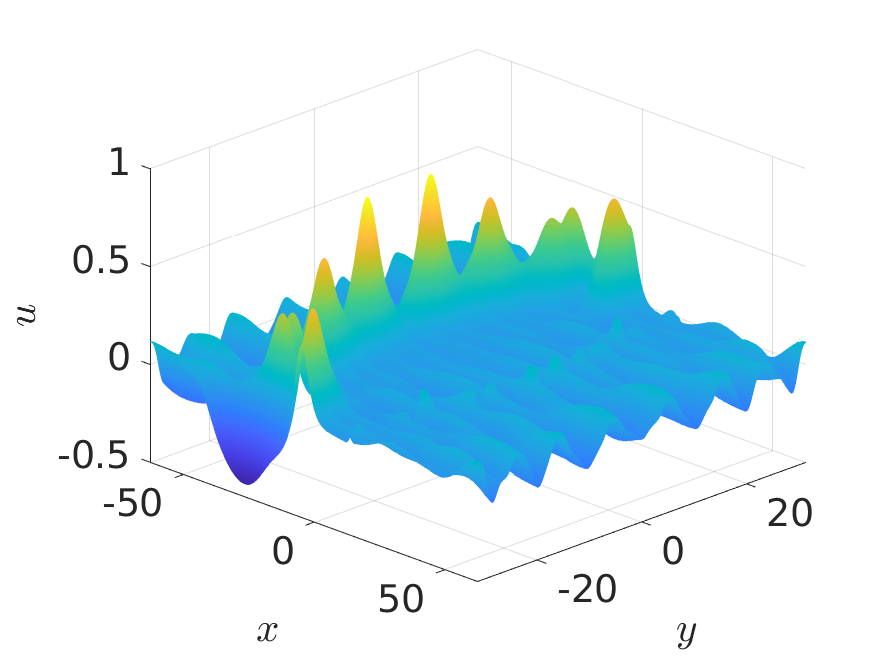}
\caption{Solution to the subcritical $(p =2)$ ZK for locally 
perturbed initial data 
(\ref{initial_gauss}) $u_0 = Q_c + a\exp^{-b^2r^2}$ 
, with $a = 0.1$, $b = 3$ and supercritical speed $c =1$:
on the left  the evolution of the $|u|_\infty$ showing stabilisation and on the right the local difference between the 
leading lump solitons and fitted lumps for $t=90$. }
\label{fig:sub_unstable_local}
\end{figure}

In this subsection we consider the unstable regime. To do this we 
study similar initial conditions as in the stable case. However, we 
change the domain and consider $L_x = L_y = 10$ with 
$N_{x}=N_{y}=2^{11}$ Fourier modes. 
We first study the localized perturbation (\ref{initial_gauss}) with $a = 0.1$ and $b = 3$
and soliton speed $c=1$ for $t\in[0,100]$ using $N_{t}=10^5$ time steps. The initial data are similar to the ones on the left of 
Fig.~\ref{fig:sub_stable_local_initial}. As $c> 
c^*=4/(5L_y^2)$ (we use $c\gg c^{*}$ in order to see lump formation on 
the chosen computational domain before radiation re-entering the 
domain having an effect on the results), the instability is clearly  visible in the $L^{\infty}$ 
norm of the solution on the left of Fig.~\ref{fig:sub_unstable_local}. 

\begin{figure}[!htb]
\includegraphics[width=0.49\hsize]{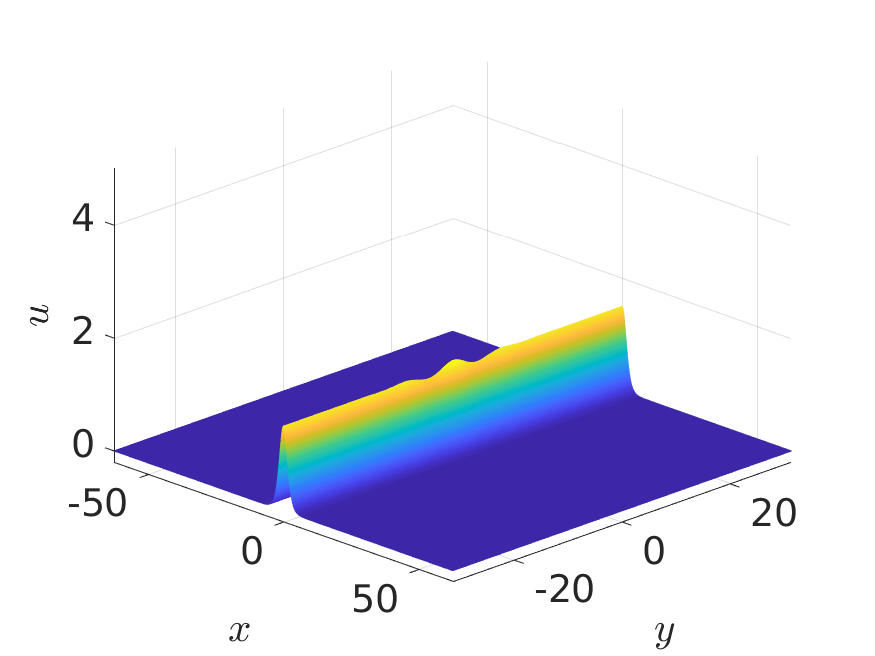} 
\includegraphics[width=0.49\hsize]{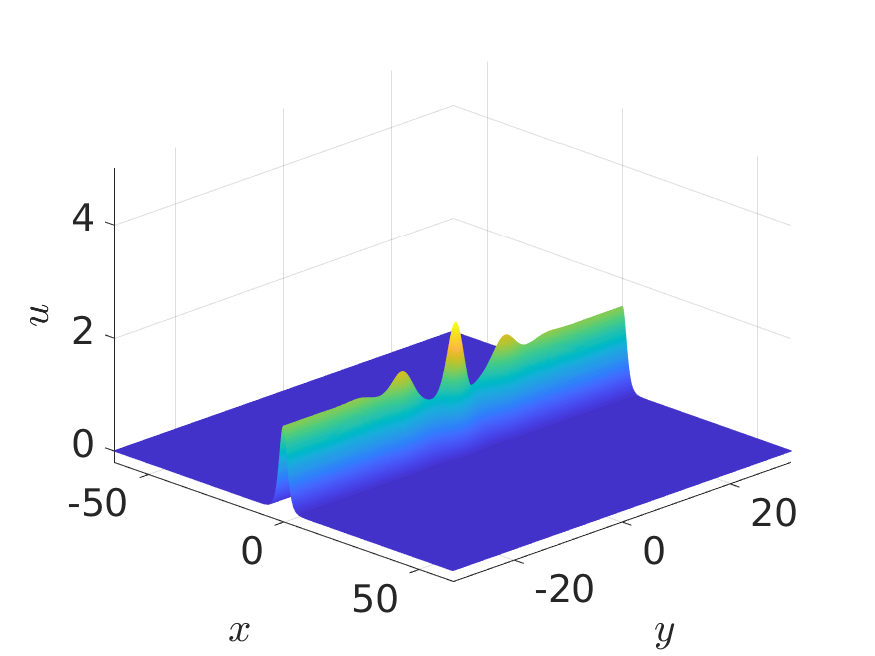}\\ 
\includegraphics[width=0.49\hsize]{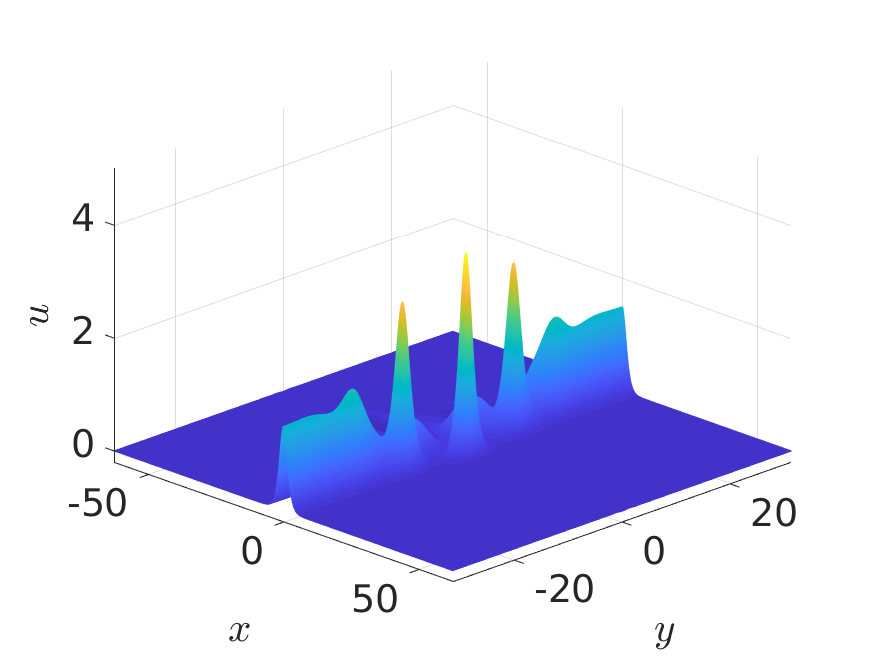} 
\includegraphics[width=0.49\hsize]{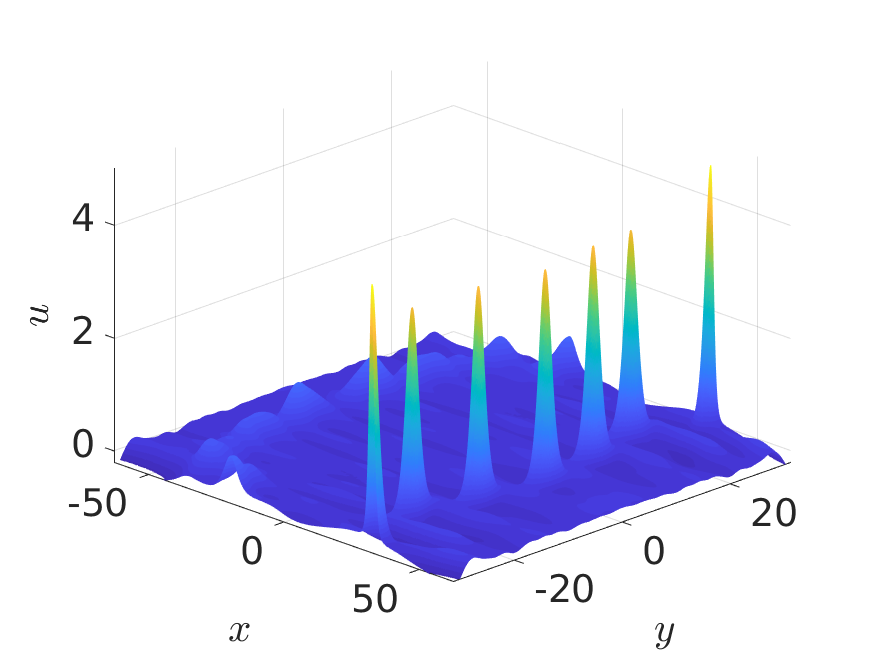} 
\caption{Subcritical ZK $(p=2)$: Evolution of the solution for locally perturbed line soliton initial data (\ref{initial_gauss}) $u_0 = Q_c + a\exp^{-b^2r^2}$ 
, with $a = 0.1$, $b = 3$ and supercritical speed $c =1$ for $t=20$, $t = 30$, $t=50$ and $t = 100$ . Note that 
the largest soliton appears on the opposite side of the initial 
perturbation and thus appears on both sides of the computational 
domain since the underlying torus is cut along this line. Note also that the orientation of the plot is changed, to provide better view to the area in front of the line soliton. }
\label{UnstableLineLoc}
\end{figure}

The initial data appear to decompose into an array of lumps as in the 
KP case, see \cite{KS12}, as can be recognized in Fig.~\ref{UnstableLineLoc}. To show that the 
array of peaks can indeed be interpreted as an array of lumps, we fit 
the solution to a multi-lump solution. As the lumps are exponentially 
localised, such a solution can be constructed as a sum of single 
soliton solutions, provided the peaks are sufficiently far away from 
each other, which is the case here (we simply subtract the solution 
$\mathcal{Q}$ of (\ref{Q}) after applying the scaling (\ref{Qscal}) 
at the maximum). The difference between the 
solution and the fitted multi-soliton is shown in 
Fig.~\ref{fig:sub_unstable_local} on the right. Locally, the difference is 
again of  the order of magnitude of the radiation, however we can see 
a possible development of a second wave of smaller solitons at the back. Mass re-entering the domain in the $x$ direction does not allow us to investigate this development further.

\begin{figure}[!htb]
 \includegraphics[width=0.45\hsize]{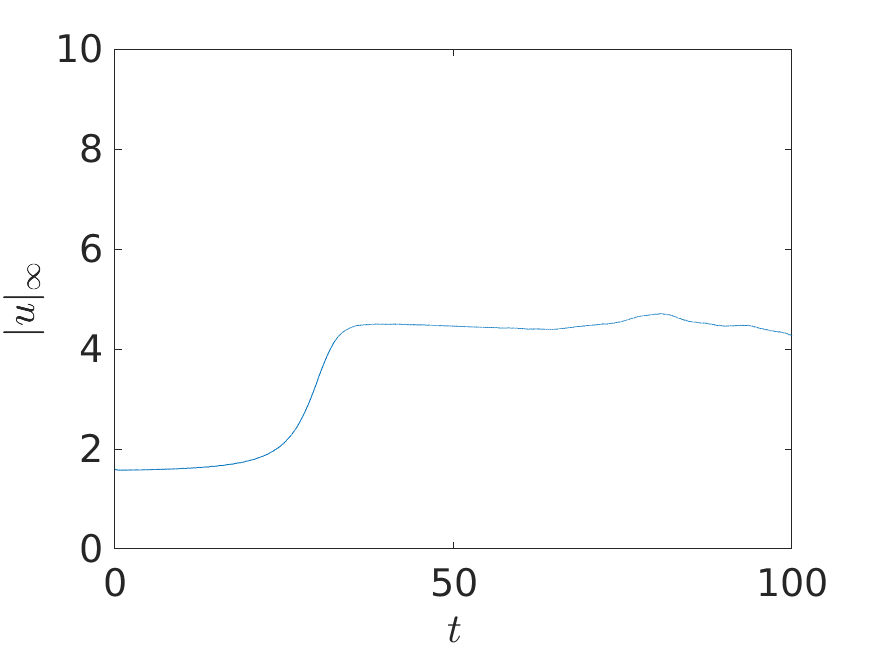} 
 \includegraphics[width=0.45\hsize]{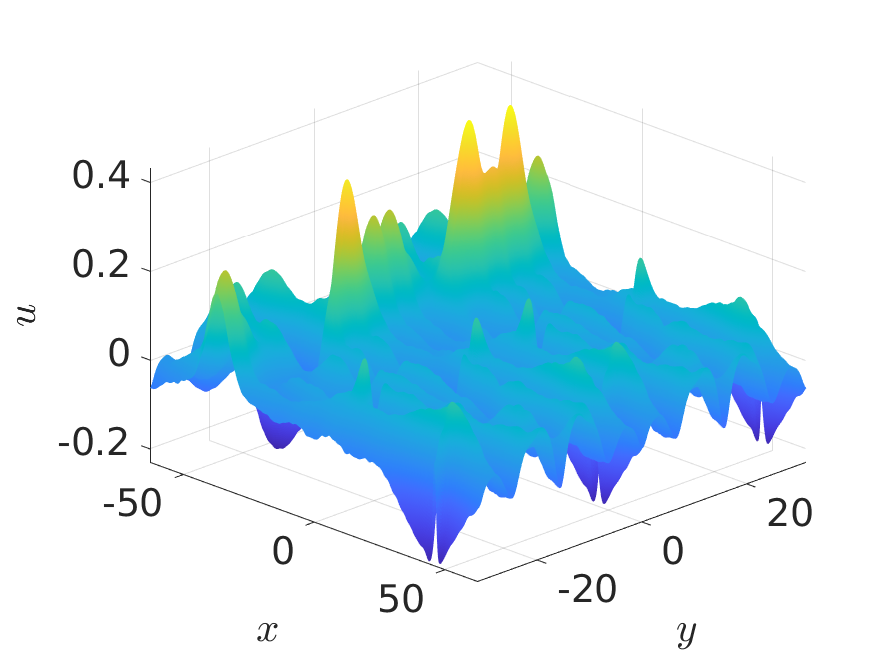}
\caption{Solution to the subcritical $(p =2)$ ZK for periodically perturbed initial data 
(\ref{initial_periodic}) $u_0 = Q_c + a\cos^2(y/L_y)\exp^{-b^2x^2}$ 
, with $a = 0.2$, $b = 1$  and supercritical speed $c=1$: 
on the left the evolution of the 
$|u|_\infty$ showing stabilisation of the solution at time $t=90$, on 
the right the local difference 
between the peaks and fitted lumps for $t=90$. }
\label{fig:Sub_unstable_periodic}

\end{figure}

Next, we study periodically perturbed initial data of the form 
(\ref{initial_periodic}) with $a =0.2$, $b = L_y$, $\delta=0$, and the same numerical parameters as in the previous case, see the left of 
Fig.~\ref{fig:sub_stable_periodic_initial}. The $L^{\infty}$ norm on 
the left of Fig.~\ref{fig:Sub_unstable_periodic} indicates the 
instability of the line soliton.  
\begin{figure}[!htb]
\includegraphics[width=0.49\hsize]{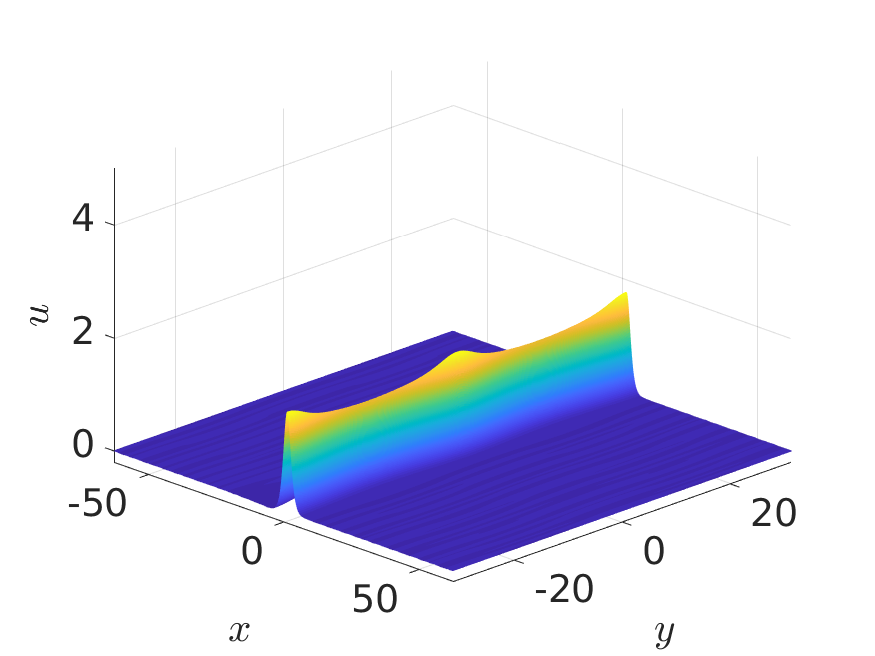} 
\includegraphics[width=0.49\hsize]{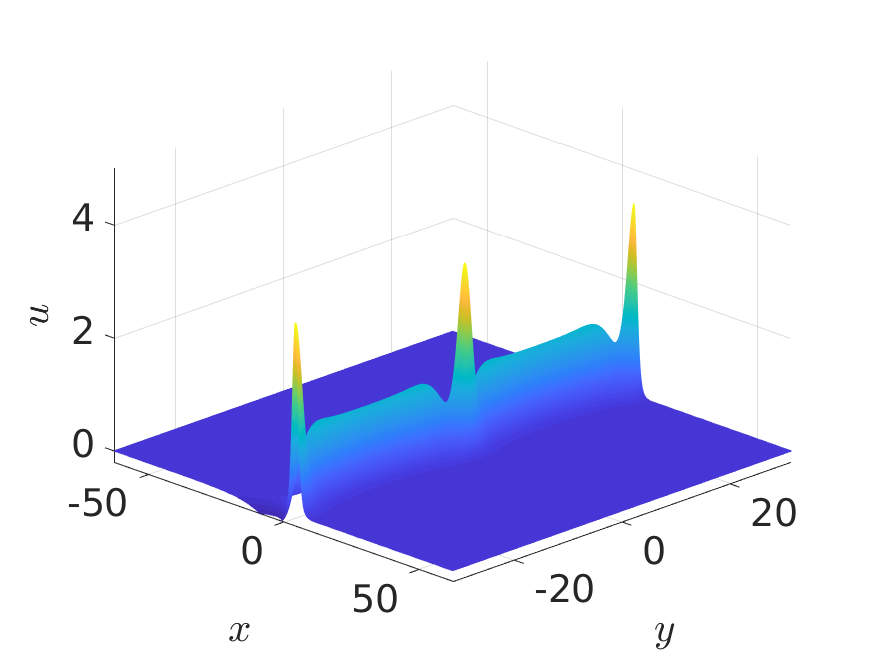}\\ 
\includegraphics[width=0.49\hsize]{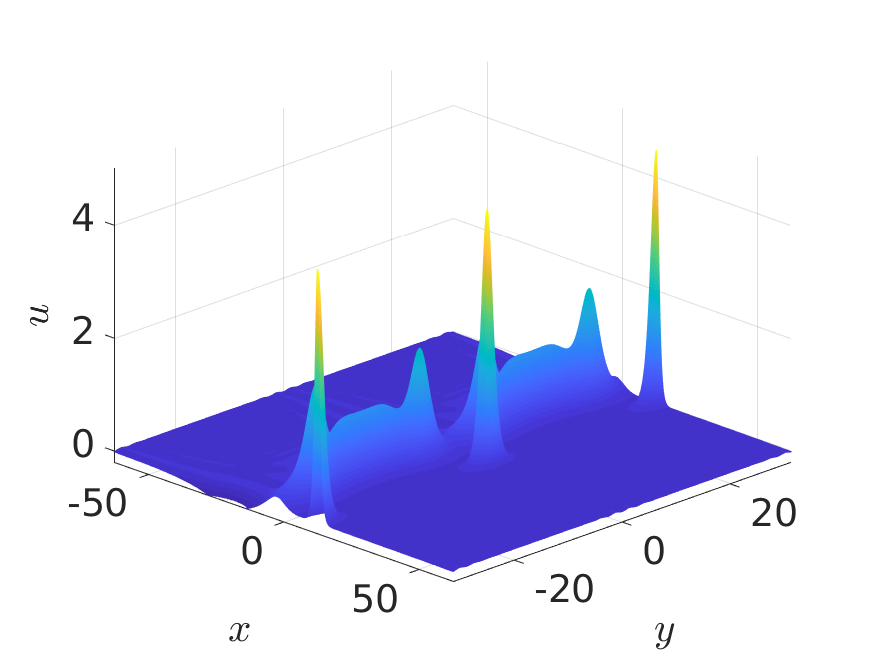} 
\includegraphics[width=0.49\hsize]{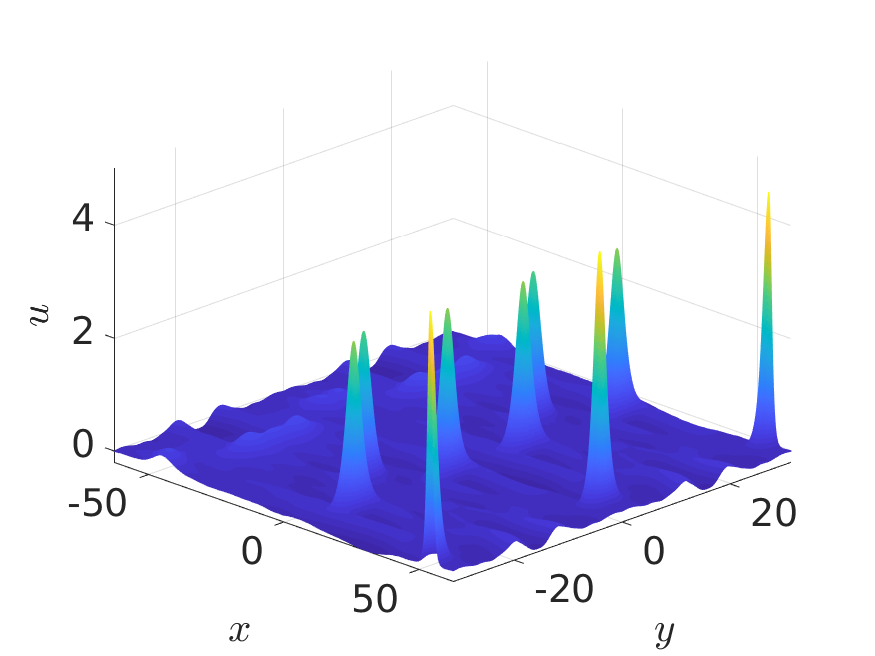}  
\caption{Subcritical ZK $(p=2)$: Evolution of the solution for periodically perturbed line soliton initial data (\ref{initial_periodic}) $u_0 = Q_c + a\cos^2(y/L_y)\exp^{-b^2x^2}$ 
, with $a = 0.2$, $b = 1$ and supercritical speed $c =1$ for $t=20$, $t = 30$, $t=50$ and $t = 90$. }
\label{UnstableLinePer}
\end{figure}

The solution for $t=90$ for these initial data is shown  
in Fig.~\ref{UnstableLinePer}. Once more an array of lumps appears 
to form. By fitting the peaks to a lump,  we can identify 
eight lump solitons, with a fit correct to 2\%, i.e., of the error of 
the radiation, as can be seen in the right of Fig.~\ref{fig:Sub_unstable_periodic}. 

\section{Critical nonlinearity}
In this section we address the $L^{2}$ critical case. 
Three variables determine the stability of the line soliton, as in 
the theorem in \cite{YAM17}: the size of the cylinder (in our case 
this is a torus, however the longitudal dimension is sufficiently 
large to treat it as a cylinder), the speed of the line soliton and 
the size of the deformation. We are not aware of analytical results in 
this direction. We consider cases where a small perturbation, that 
is, a perturbation with a norm  within 10\% of the infinity norm of 
the line soliton.
We expect two regimes as before though the critical speed $c^{*}$, below 
which the line soliton is stable, is not 
known in this case. Numerically it is difficult to determine the critical speed, however we show that a sufficiently small and slow line 
soliton is stable. For larger values of $c$, lump formation is
possible, however since the lumps are conjectured to be unstable against 
blow-up, see \cite{KRS}, an $L^{\infty}$ blow-up as in 
Conjecture \ref{C:2} of \cite{KRS} can 
be also observed here. We do not go into details here since the 
results appear to be in accordance with what has been found in 
\cite{KRS}. 

\begin{figure}[!htb]
\includegraphics[width=0.45\hsize]{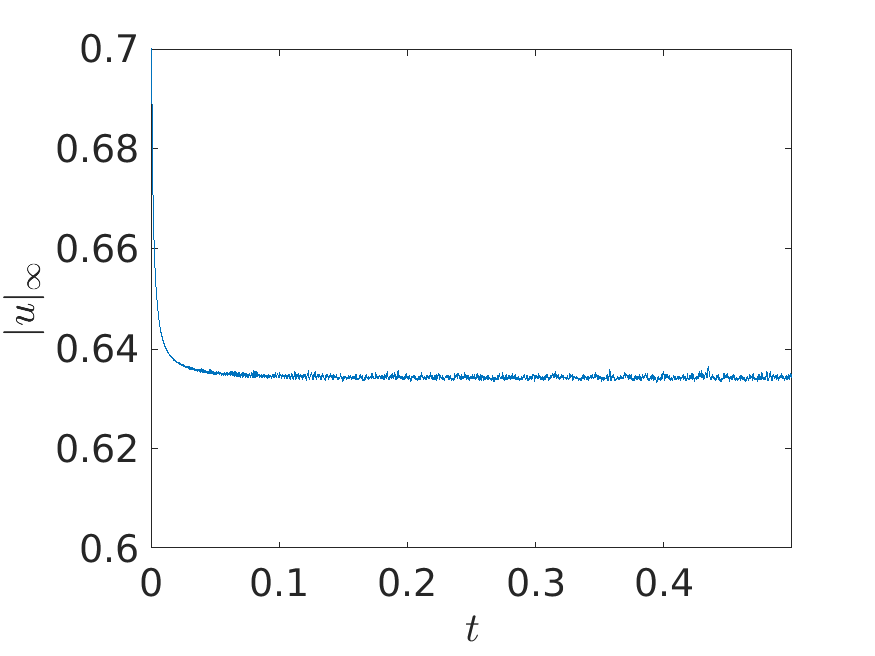}
\includegraphics[width=0.45\hsize]{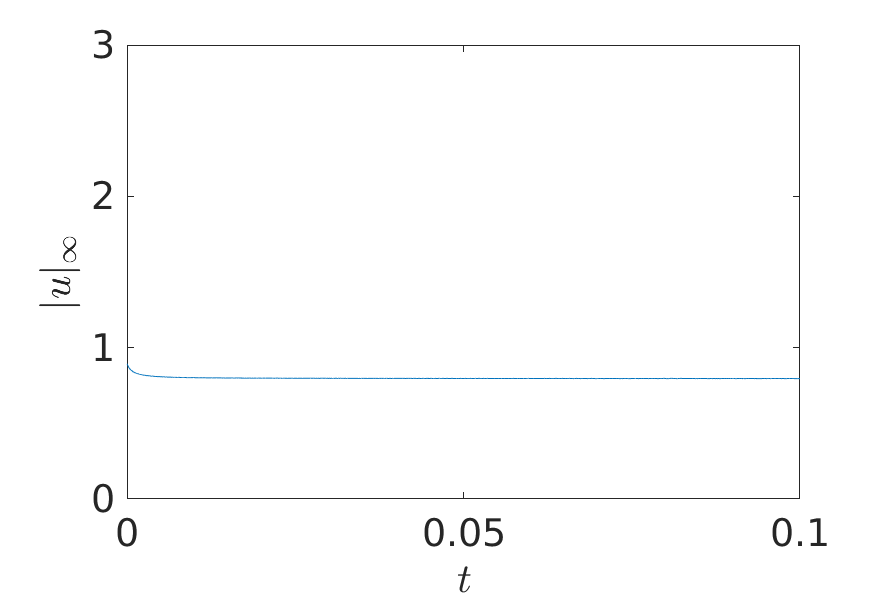}
 \caption{$L^{\infty}$ norm of the solution for the locally perturbed line soliton initial data 
(\ref{initial_gauss}) $u_0 = Q_c + a\exp^{-b^2r^2}$, with $a = 0.1$, $b = 3$, and subcritical speed 
$c = 0.2$: on the left
for the critical $(p =3)$ ZK equation, on the right for the super-critical nonlinearity $(p=4)$. A domain with $L_y = 1$ is sufficient to capture all the dynamics.  }
\label{fig:subcritinf}
\end{figure}

If the speed $c$ of the line soliton is smaller than some unknown 
critical value $c^{*}$, we expect the line soliton to be stable under 
small perturbations. This is 
indeed shown by a numerical experiment  on the left of Fig.  
\ref{fig:subcritinf}. We see for the $L^{\infty}$ norm that a small line soliton with 
$c=0.2$  stabilises under a small local Gaussian perturbation of the 
form (\ref{initial_gauss}). Since 
the latter is finite in order to see numerical effects for finite 
computing times, the solution tends to a  line 
soliton with slightly different speed. The difference in height between the unperturbed initial 
soliton and the final one is about 2\%.

\begin{figure}[!htb]
\includegraphics[width=0.49\hsize]{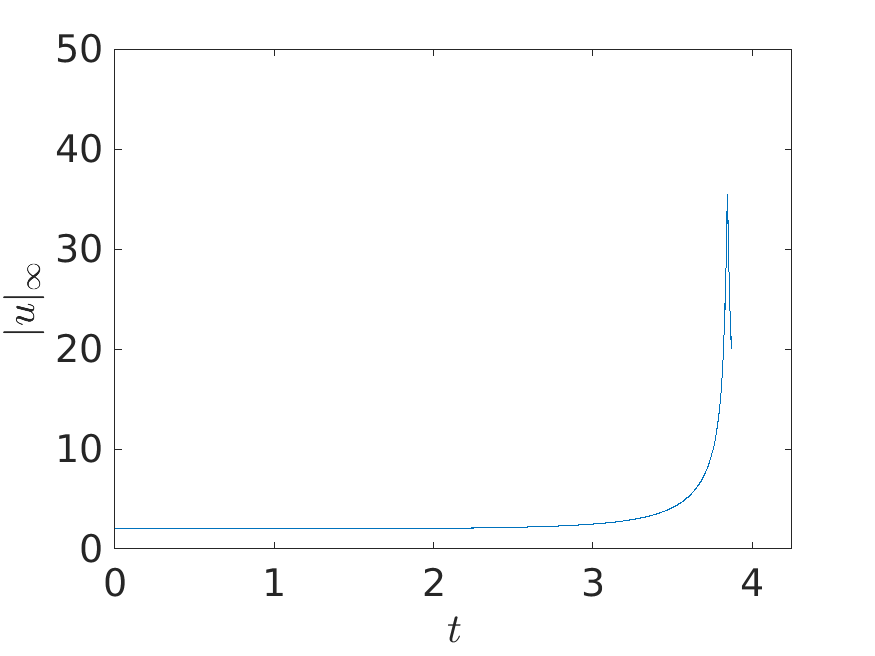} 
 \includegraphics[width=0.49\hsize]{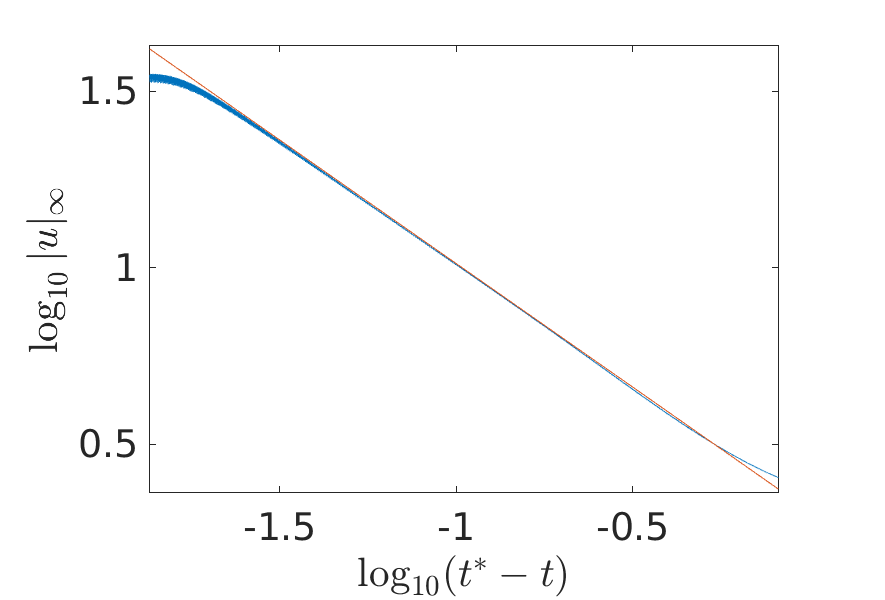}
\caption{Solution to the critical $(p =3)$ ZK equation for locally deformed soliton initial data 
(\ref{initial_gauss}) $u_0 = Q_c + a\exp^{-b^2r^2}$, with $a = 0.1$, $b = 3$,  and $c=1$:  on the  left 
the evolution of the $L^{\infty}$ norm, on the right a fit of the 
$L^{\infty}$ norm for approximately $4*10^4$ time steps 
according to (\ref{fit}) to obtain $t^* = 3.858$ and blow-up rate $q = 0.7$.}
\label{fig:crit_stable_local}
\end{figure}

\begin{figure}[!htb]
\includegraphics[width=0.49\hsize]{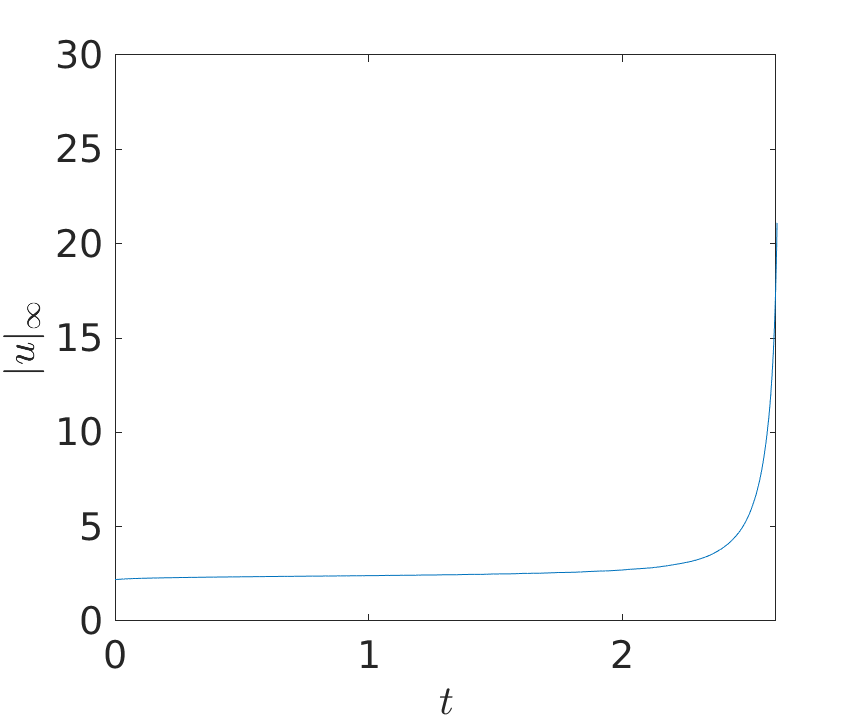} 
 \includegraphics[width=0.49\hsize]{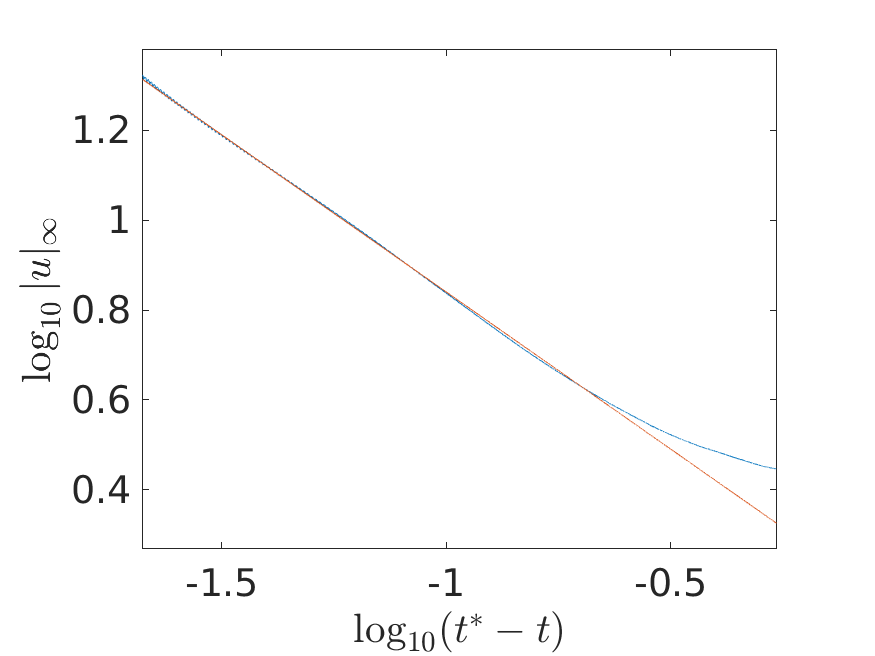}
\caption{Solution to the critical $(p =3)$ ZK equation for locally deformed soliton initial data 
(\ref{initial_gauss}) $u_0 = Q_c + a\exp^{-b^2r^2}$, with $a = 0.1$, $b = 3$ and superciritical speed $c =2$: on the left  the evolution of $|u|_{\infty}$, on the 
right a fit of the $L^{\infty}$ norm for the last approximately 
$1.5*10^4$  time steps to find $t^* = 2.631$ and blow-up rate $q = 0.7$ 
according to (\ref{fit}).}
\label{CriticalLineLoc_fig}
\end{figure}

In contrast, if the speed 
 $c>c^{*}$ we expect the line soliton to be unstable against lump 
 formation that are in turn unstable against a blow-up.  As the numerical 
 code cannot come arbitrarily close to a singularity without losing precision, we can track 
the solution only until close to the time of formation of a first singularity. Whether or 
not the solution forms additional singularities at the location of 
other lumps at a later time cannot be answered with the present means. We will only 
address below the question of the first singularity to appear. 
 We are able to recover the blow-up rate by using a multi parameter fit on the evolution of 
$|u|_\infty$  according to (\ref{resc}) to the 
expression,              
\begin{equation}
\log_{10}|u|_\infty = q\log_{10}(t^* - t) + r,
\label{fit}
\end{equation}      
where we fit for the critical time $t^*$, blow-up rate $q$ and 
parameter $r$ using least squares. We use data points that are sufficiently close to the 
blow-up, but for which the conservation of mass is still sufficient. The fitting parameters given in the caption of Fig.~\ref{fig:crit_stable_local} 
do not depend significantly on the number of points used for the 
fitting. In our plots on Fig. \ref{fig:crit_stable_local} and 
\ref{CriticalLineLoc_fig} on the right we include points that are not 
used for the fitting, in order to better illustrate the three 
regimes: non-critical regime on the right of the figure, when the solution is far away from the blow-up point, followed by a section that is sufficiently close to the blow-up, and the solution is well tracked and the red line and the blue curve match, and finally on the left of the plot the numerical scheme loses the solution.  
\begin{remark}
	The fitting (\ref{fit}) only gives approximate values of the 
	blow-up rate $q$ even in situation for which theorems exist as for 
	gKdV, see for instance \cite{KP2}, since blow-up time $t^{*}$ and 
	rate $q$ have to be determined in one fit. The found values have 
	to be seen as indicative.
\end{remark}

\begin{figure}[!htb]
\includegraphics[width=0.49\hsize]{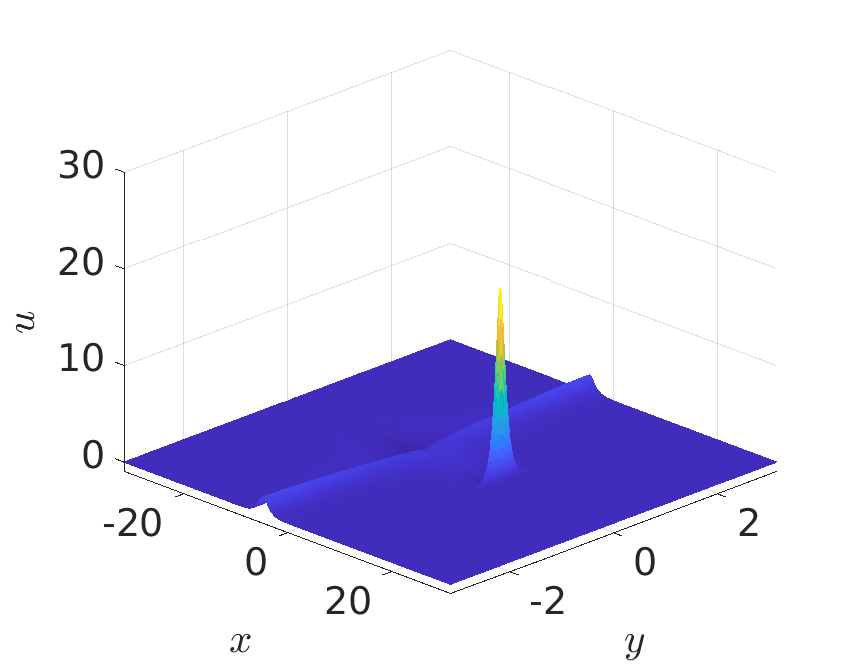} 
 \includegraphics[width=0.49\hsize]{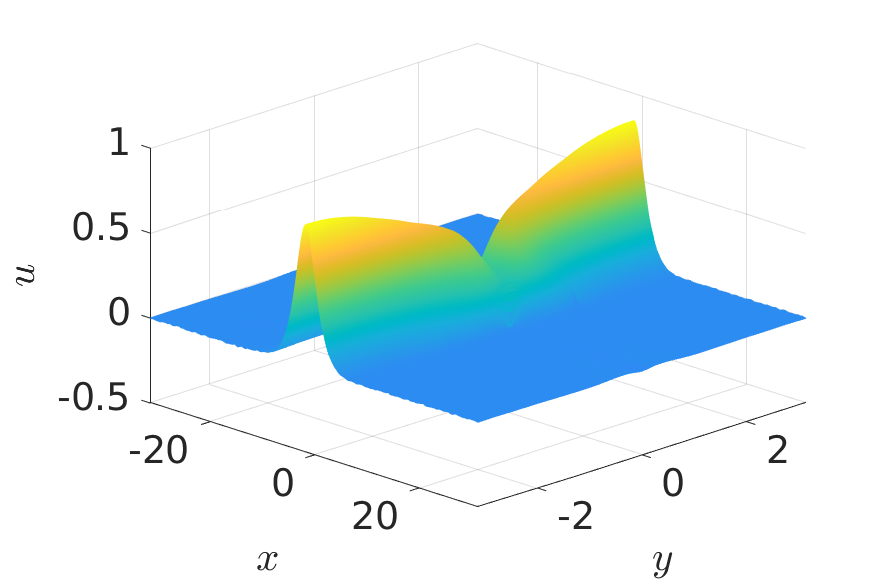}
\caption{Solution to the critical $(p =3)$ ZK equation for locally deformed soliton initial data 
(\ref{initial_gauss})$u_0 = Q_c + a\exp^{-b^2r^2}$, with $a = 0.1$, $b = 3$ and supercritical speed $c = 2$ with $L_y=1$:   on the left the solution for time 
$t=3.75$ and the local difference between the peak and a dynamically 
rescaled lump on the right. }
\label{fig:critical_local_plots}
\end{figure}

The code is stopped once the numerically computed mass is no longer 
conserved relatively to better than $10^{-3}$. The solution for 
$t=3.75$ is shown on the left of Fig.~\ref{fig:critical_local_plots}. On the 
right of the same figure, we show the difference between the peak and a 
dynamically rescaled lump. It can be seen that the lump gives 
the blow-up profile to the order of 2\% though the singularity is 
not yet reached. Note that the singularity is formed at the same $y$ 
position as the initial (positive) deformation. A negative 
deformation of sufficient mass produces a blow-up in the opposite $y$ 
- direction. 

We now look at a periodic perturbation on a large torus with period 
$2\pi L_{y}$ in $y$-direction. Here we take $L_y =10$, so that the expected two lumps are not in a strong interaction regime, see \cite{KRS} for fusion of ZK lump solutions, and initial condition 
(\ref{initial_periodic})

where $a =0.1$, $b =L_y$, $c = 2$ and $\delta=\pi/2$.
The evolution of the $L^{\infty}$ norm of the solution for these 
initial data as well as the blow-up rate are shown on  Fig.~\ref{CriticalLineLoc_fig}. We note that the solution blows up simultaneously at two spatial points and that the $y$ coordinates of these two points are where the initial deformation is minimal.  

\begin{figure}[!htb]
\includegraphics[width=0.49\hsize]{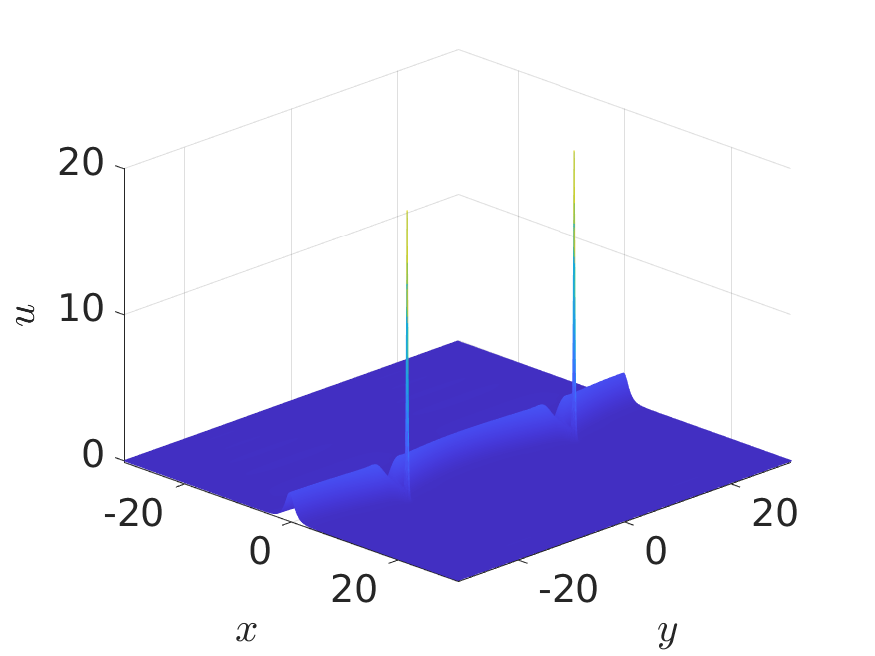} 
 \includegraphics[width=0.49\hsize]{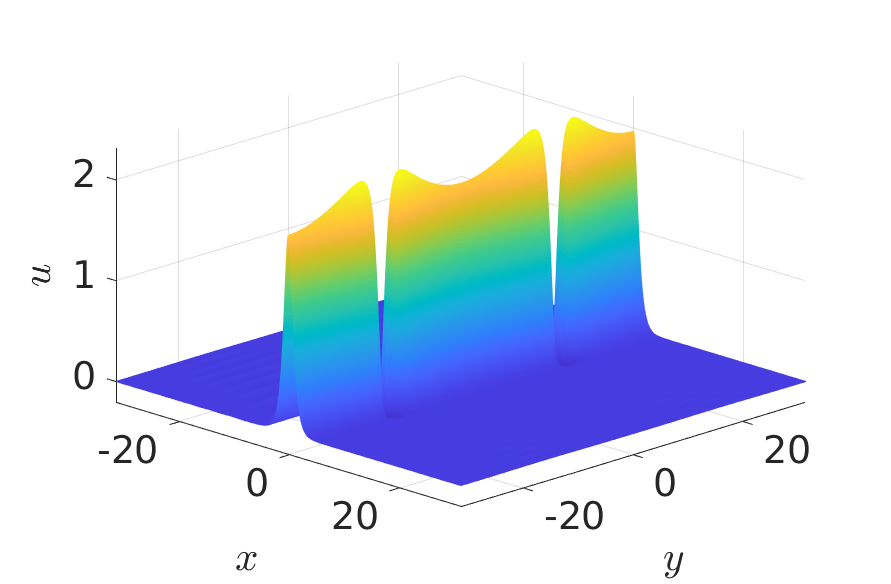}
\caption{Solution to the critical $(p =3)$ ZK equation for locally deformed soliton initial data 
(\ref{initial_periodic}) $u_0 = Q_c + a\cos^2(y/L_y+ \delta)\exp^{-b^2x^2}$, $\delta = \pi/2$ and critical speed $c=2$:  on the left the solution $u$ at time 
$t=2.6$, on the right the local difference between the peak and a 
fitted lump. }         
\label{CriticalLineLoc_fig2}
\end{figure}

The code is stopped once a lack in mass conservation indicates a loss 
of accuracy. We show on the left of Fig.~\ref{CriticalLineLoc_fig2} 
the solution at this time, $t=2.6$. On the right of the same figure 
we show the difference between this solution and lumps fitted to the 
peaks. 

\section{Super-critical nonlinearity}
In the super-critical case we expect a similar behaviour as in the 
critical case, namely two regimes depending on the speed of the 
initial line soliton. If $c<c^{*}$, where $c^{*}$ is once more 
unknown, the final state is a line soliton of slighly 
different mass. This is indeed shown by a numerical experiment 
described in Fig.~\ref{fig:subcritinf} on the right, that tracks the 
$L^{\infty}$ norm for the initial data (\ref{initial_gauss}) with 
$a=0.1$ and $b=3$
with $c = 0.2$.
 We see that a small 
line soliton perturbed by a small Gaussian leads to a final state 
being a line soliton of slightly different speed. The difference in 
height between the unperturbed initial condition and the final is about 2\%. We observe the 
same behaviour for different types of deformations e.g. localised in 
$x$ and periodic in $y$ or localised negative deformations.      

\begin{figure}[!htb]
\includegraphics[width=0.45\hsize]{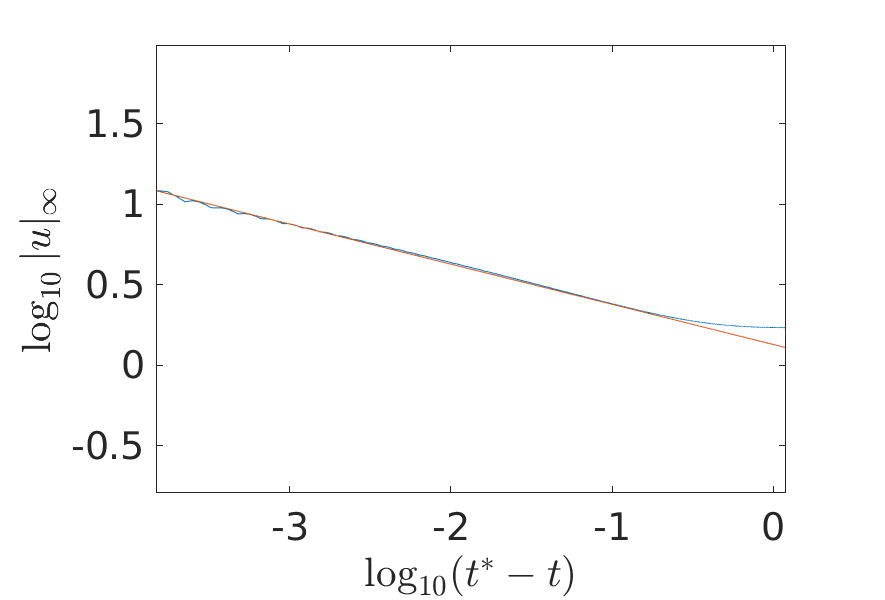} 
 \includegraphics[width=0.45\hsize]{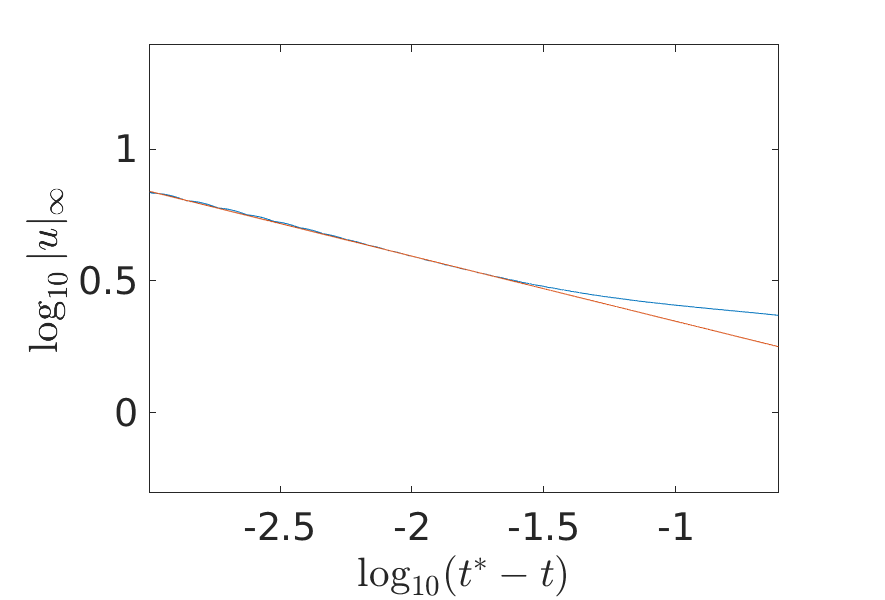}
\caption{Solution to the supercritical $(p =4)$ ZK equation: Blow-up 
rate   for locally perturbed (left) initial data $u_0 = Q_c + a\exp^{-b^2r^2}$ with $a= 0.1$, $b = 3$ and $c = 2$ showing blow-up rate 
$q = 0.249$ and critical time $t^* = 1.2134$,  and periodically 
perturbed data (\ref{initial_periodic}) $u_0 = Q_c + a\cos^2(y/L_y +\gamma)\exp^{-b^2x^2}$ with $a= 0.1$, $b = 3$, $\gamma = \pi/2$, and $c = 2$ (right), exhibiting simultaneous blow-up in two points with blow-up rate $q = 0.247$ and critical time $t = 0.611$. }
\label{SuperLineLoc_fig}
\end{figure}

For supercritical speeds, the behaviour is similar to the critical 
case, and again we observe simultaneous blow-up at two different 
spatial points in the periodically deformed case for periodically 
perturbed initial data (\ref{initial_periodic}), see 
Fig.~\ref{SuperLineLoc_fig} on the right. The behavior for a Gaussian 
perturbation (\ref{initial_gauss}) is shown 
on the left of the same figure. 
The observed 
blow-up rates are close to $1/3$ as predicted by the self-similar 
mechanism conjectured in \cite{KRS}. 

\begin{figure}[!htb]
\includegraphics[width=0.49\hsize]{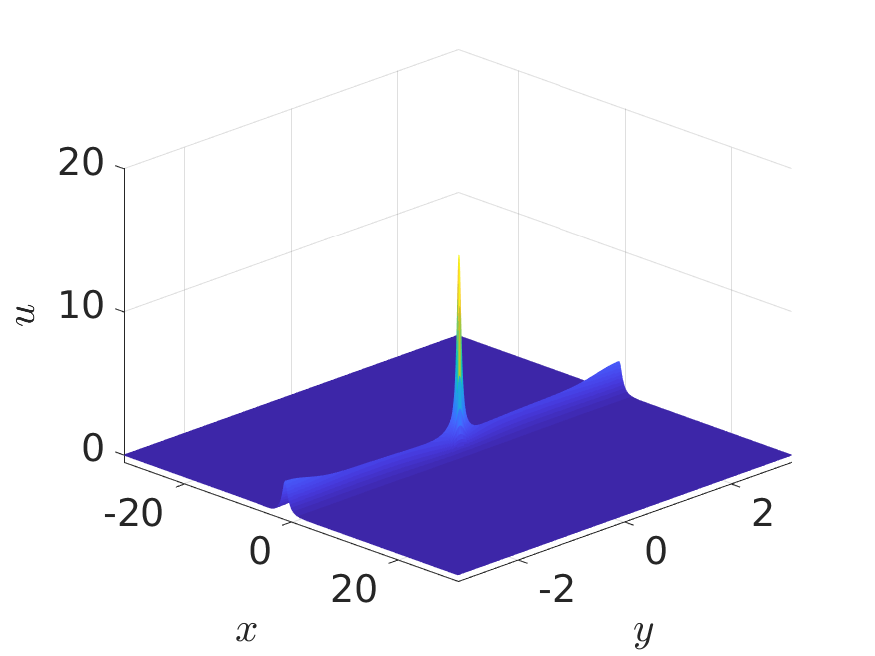}
\includegraphics[width=0.49\hsize]{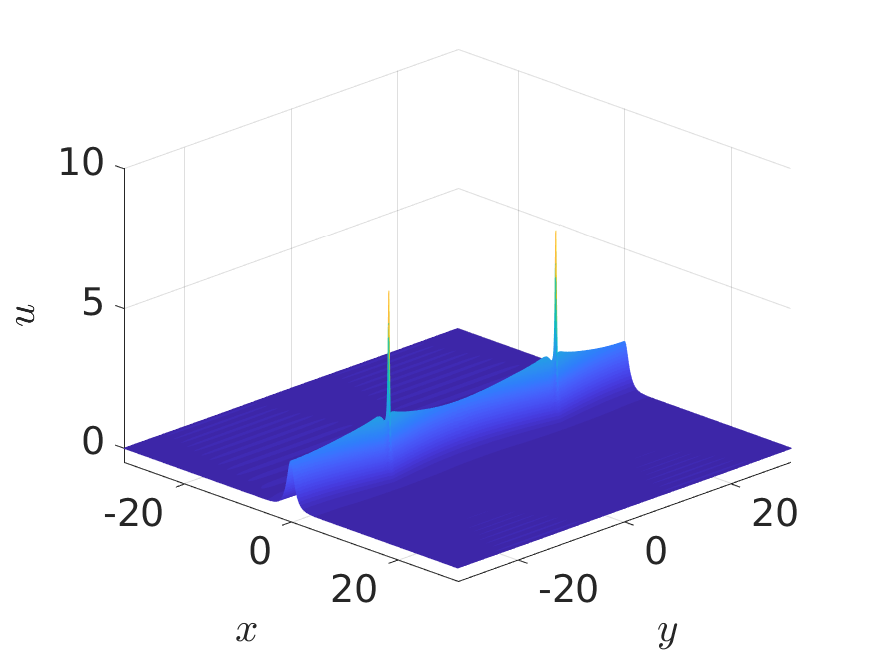}
\caption{Profile close to the blow-up for the solution to the supercritical ZK equation $(p=4)$ for locally perturbed (left) soliton initial data 
(\ref{initial_gauss})$u_0 = Q_c + a\exp^{-b^2r^2}$ with $a= 0.1$, $b 
= 3$ and $c = 2$  at $t = 1.2126$, and periodically perturbed initial data (\ref{initial_periodic}), $u_0 = Q_c + a\cos^2(y/L_y +\gamma)\exp^{-b^2x^2}$ with $a= 0.1$, $b = 3$, $\gamma = \pi/2$ (right) at time $t = 0.61$. }
\label{fig:SuperProfile}
\end{figure}

We show the solutions close to blow-up in 
Fig.~\ref{fig:SuperProfile}, on the left for a localized 
perturbation, on the right for a deformed line soliton. It can be 
seen that lump-like structures form in both cases which will eventually blow up. 

\section{Conclusion}
In this paper we have presented a detailed numerical study of 
perturbations of line solitons on $\mathbb{T}^{2}$ approximating a 
situation on $\mathbb{R}\times \mathbb{T}_{L}$. We consider 
perturbations being localised both in $x$ and $y$ and perturbations 
localised in $x$, but periodic in $y$, i.e., deformed line solitons. It was shown that 
the line soliton is stable for $c$ smaller than some critical 
speed $c^{*}$, also in $L^{2}$ critical and supercritical cases (but 
subcritical in 1D). 
In the latter cases, the precise value of $c^{*}$ appears to be 
unknown. 
For values of $c>c^{*}$, the line soliton appears to be unstable 
against lump formation. Since the lumps are known to be strongly 
unstable in the critical and supercritical cases, a blow-up in finite 
time was 
observed in these settings. 

An interesting question would be to determine the value of the 
critical speed $c^{*}$ which is numerically difficult, but could 
be accessible analytically. Since the ZK equation can also be studied 
in 3D, the 2D lumps can be extended to line solitons of the 3D ZK 
equation. The transverse stability of these solutions will be studied 
numerically in an ensuing work.


\begin{thebibliography}{9}
\bibitem{AKS} 
J. Arbunich, C. Klein, C. Sparber, 
\emph{On a class of derivative Nonlinear Schr\"odinger-type equations in two spatial dimensions},  
M2AN 53(5), (2019), 1477 - 1505.
\bibitem{BL}
H. Berestycki and P.-L. Lions,
\emph{Nonlinear scalar field equations},
Arch. Rational Mech. Anal. {\bf 82} (1983), 313-376.

\bibitem{CMPS} 
R. C\^{o}te, C. Mu\~{n}oz, D. Pilod, and G. Simpson, 
\emph{Asymptotic stability of high-dimensional Zakharov-Kuznetsov solitons}, 
Arch. Ration. Mech. Anal. 220 (2016), no. 2, 639--710.

\bibitem{CM}
S. Cox and P. Matthews, 
\emph{Exponential Time Differencing for stiff Systems}, 
J. of Comp. Phys., 176 (2002), 430-455.

\bibitem{dB} 
A. de Bouard,  
\emph{Stability and instability of some nonlinear dispersive solitary waves in higher dimension}, Proc. Roy. Soc. Edinburgh Sect. A 126 (1996), no. 1, pp. 89--112.

\bibitem{Bridges}
T.J. Bridges,
\emph{Universal geometric conditions for the transverse instability of solitary waves}, Phys. Rev. Lett. {\bf 84} (2000), 2614-2617.

\bibitem{F95}  
A. V. Faminskii, 
\emph{The Cauchy problem for the Zakharov-Kuznetsov equation}. (Russian) 
Differentsialnye Uravneniya 31 (1995), no. 6, 1070--1081, 1103; translation in Differential Equations 31 (1995), no. 6, 1002--1012.

\bibitem{FLP} 
L.G. Farah, F. Linares and A. Pastor, 
\emph{A note on the 2D generalized Zakharov-Kuznetsov equation: Local, global, and scattering results}, J. Diff. Eq. 253 (2012), 2558--2571.

\bibitem{FHR2} 
L. G. Farah, J. Holmer and S. Roudenko, 
\emph{Instability of solitons - revisited, II: the supercritical Zakharov-Kuznetsov equation}, 
Contemp. Math., 725, Amer. Math. Soc., 89–109.

\bibitem{FHR3} 
L. G. Farah, J. Holmer and S. Roudenko, 
\emph{Instability of solitons in the 2d cubic Zakharov-Kuznetsov equation}, 
Fields Institute Communications, vol 83 (2019), 
Eds: Miller P., Perry P., Saut JC., Sulem C., Nonlinear Dispersive Partial Differential Equations and Inverse Scattering.  Springer, New York, NY

\bibitem{FHRY} 
L. G. Farah, J. Holmer, S. Roudenko and Kai Yang, 
\emph{Blow-up in finite or infinite time of the 2D cubic Zakharov-Kuznetsov equation}, 
arXiv:1810.05121

\bibitem{Han-Kwan}
D. Han-Kwan,
\emph{From Vlasov-Poisson to Korteweg and Zakharov-Rubenchik},
Comm. Math. Phys. {\bf 324}, no.3 (2013), 961-993. 

\bibitem{HO}
M. Hochbruck, A. Ostermann, 
\emph{Exponential integrators}, 
Acta Numerica (2010), pp. 209-286, doi:10.1017/S0962492910000048

\bibitem{ITK}
H. Iwasaki, S. Toh and T. Kawahara,
\emph{Cylindrical quasi-solitons of the Zakharov-Kuznetsov equation},
Physica D {\bf 43} (1990), 293-303.

\bibitem{KK} 
A.~Kazeykina and C.~Klein, 
\emph{Numerical study of blow-up and stability of line solitons for the Novikov-Veselov equation}, 
Nonlinearity 30, 2566-2591  (2017)

\bibitem{K2021}
S. Kinoshita,
\emph{Global well-posedness for the Cauchy problem of the Zakharov-Kuznetsov equation in 2D},
Ann.I.H. Poincar\' e-AN {\bf 38} (2021), 451-505.

\bibitem{etna} 
C.~Klein,  
\emph{Fourth order time-stepping for low dispersion Korteweg-de 
Vries and nonlinear Schr\"odinger equation},  
ETNA Vol. 29 116-135 (2008).

\bibitem{KP1} 
C. Klein and R. Peter, 
\emph{Numerical study of blow-up in solutions to generalized Kadomtsev-Petviashvili equations}, Discrete Contin. Dyn. Syst. Ser. B 19 (2014), 1689-1717.

\bibitem{KP2} 
C. Klein and R. Peter, 
\emph{Numerical study of blow-up in solutions to generalized Korteweg-de Vries equations}, 
Phys. D 304 (2015), 52-78.


\bibitem{KR} 
C.~Klein and K.~Roidot,  
\emph{Fourth order time-stepping for Kadomtsev-Petviashvili and Davey-Stewartson equations},  
SIAM J. Sci. Comput., 33(6), 3333-3356. DOI: 10.1137/100816663 (2011). 

     \bibitem{KRS} C. Klein, S. Roudenko, N. Stoilov, Numerical study of 
     Zakharov-Kuznetsov equations in two dimensions,  J. Nonl. Sci  
	 31(26) (2021) https://doi.org/10.1007/s00332-021-09680-x
     

	 \bibitem{KRS2} C. Klein, S. Roudenko, N. Stoilov, Numerical study of 
	 soliton stability, resolution and interactions in the 3D 
	 Zakharov-Kuznetsov equation, Physica D 423 (2021) 132913,
	 https://doi.org/10.1016/j.physd.2021.132913

\bibitem{KS12} C.~Klein and J.-C.~Saut,  \emph{Numerical study of blow up and 
stability of solutions of generalized Kadomtsev-Petviashvili 
equations}, J. Nonl. Sci. Vol. 22 (5), 763-811 (2012). 

	 
\bibitem{KS}  
C.~Klein and N.~Stoilov, 
\emph{A numerical study of blow-up mechanisms for Davey-Stewartson II systems}, 
Stud. Appl. Math., DOI : 10.1111/sapm.12214 (2018)

\bibitem{Kuz}E.A. Kuznetsov, Stability criterion for solitons of the 
Zakharov–Kuznetsov-type equations, Physics Letters A 382 (2018) 2049–2051

\bibitem{KRZ}E.A. Kuznetsov, A.M. Rubenchik, V.E. Zakharov, Soliton 
Stability in plasmas and hydrodynamics, Phys. Rep. 142(3) (1986) 
103-165. 

\bibitem{Kw}
M.K. Kwong,
\emph{Uniqueness of positiv eradial solutions of $\Delta u-u+u^p=0$ in $\R^n$,
Arch. Rational Mech. Anal. {\bf 105} (1989), 243-266.}



\bibitem{Lan}Yang Lan, Stable Self-Similar Blow-Up Dynamics for 
Slightly $L^{2}$-Supercritical Generalized KDV Equations, 
Communications in Mathematical Physics volume 345, pages 223–269 
(2016)  

\bibitem{LLS} 
D. Lannes, F. Linares and J.-C. Saut, 
\emph{The Cauchy problem for the Euler-Poisson system and derivation of the Zakharov-Kuznetsov equation}, 
Prog. Nonlinear Diff. Eq. Appl., 84 (2013), 181--213.


\bibitem{LP2009} 
F. Linares and A. Pastor, 
\emph{Well-posedness for the two-dimensional modified Zakharov-Kuznetsov equation}, 
SIAM J. Math. Anal. 41,  no. 4 (2009), 1323--1339. 

\bibitem{LPS}
F.Linares, A. Pastor and J.-C. Saut,
 \emph{Well-posedness for the ZK equation in a cylinder and on the background of a KdV soliton}, Comm. Partial Diff. Eq.  {\bf 35} 9 (2010), 1674-1689.

\bibitem{MM} 
Y. Martel and F. Merle, 
\emph{Blow up in finite time and dynamics of blow up solutions for the $L^2$-critical generalized KdV equation}, 
J. Amer. Math. Soc. 15 (2002), 617--664.

\bibitem{MM1989} 
S. Melkonian and S. A. Maslowe, 
\emph{Two dimensional amplitude evolution equations for nonlinear dispersive waves on thin films}, 
Phys. D {\bf 34} (1989), pp. 255--269.

\bibitem{M} 
F. Merle, 
\emph{Existence of blow-up solutions in the energy space for the critical generalized KdV equation}, 
J. Amer. Math. Soc. {\bf 14}, no. 3, (2001), 555--578.  

\bibitem{MoPi}
L.Molinet and D. Pilod, 
\emph{Bilinear Strichartz estimates for the Zakharov-Kuznetsov equation and applications},
 Ann.I. H. Poincar\' e-AN {\bf 32} (2015), 347-371.

\bibitem{MP1999} 
S. Munro and E. J. Parkes, 
\emph{The derivation of a modified Zakharov-Kuznetsov equation and the stability of its solutions}, 
J. Plasma Phys. {\bf 62} (3) (1999), 305--317.

\bibitem{Pel}
D. Pelinovsky
\emph{Normal form for transverse instability of the line soliton with a nearly critical speed of propagation}
Math.Model.Nat. Phenom. {\bf 13} (2018), 1-20.

\bibitem{Pu}
X. Pu,
\emph{Dispersive limit of the Euler-Poisson system in higher dimensions},
SIAM J. Math. Anal. {\bf 45} no. 2 (2013), 834-878.

\bibitem{RV2013}  
F. Ribaud and S. Vento, 
\emph{A note on the Cauchy problem for the 2D generalized Zakharov-Kuznetsov equations}, 
C. R. Math. Acad. Sci. Paris 350 (2012), no. 9-10, 499--503. 

\bibitem{RT}F. Rousset and N. Tzvetkov, Stability and instability of 
the KdV solitary wave under the KP-I flow, Comm. Math. Phys., {\bf 313} 
(2012), no. 1, 155–173.   

\bibitem{RT2} F. Rousset and N. Tzvetkov, Transverse nonlinear  instability of solitary waves for some Hamiltonian PDE's, J. Math.Pures Appl. {\bf 90} (2008), 550-590.

\bibitem{gmres} 
Y. Saad and M. Schultz, 
\emph{GMRES: a generalized minimal residual algorithm for solving nonsymmetric linear systems}, SIAM J. Sci. Comput. {\bf 7} (1986), 856-869.

\bibitem{SB2000}
R. Sipcic, D. J. Benney, 
\emph{Lump Interactions and Collapse in the Modified Zakharov-Kuznetsov equation}, Stud. Appl. Math. {\bf 105} (4) (2000), 385--403.

\bibitem{SS1999}
C. Sulem, P.-L. Sulem, 
\emph{The nonlinear Schr\"odinger equation. Self-focusing and wave-collapse.} 
Springer, 1999.

\bibitem{V2020}
F. Valet, 
\emph{Asymptotic K-soliton-like solutions of the Zakharov-Kuznetsov type equations}, 	Transactions  AMS {\bf 374} (5) (2021), 3177-3213.

\bibitem{YAM17}Y. Yamazaki, Stability for line solitary waves of Zakharov-Kuznetsov equation, Journal of Differential Equations, Volume {\bf 262}, Issue 8, (2017), 4336-4389.

\bibitem{YAM20} Y. Yamazaki, Center stable manifolds around line solitary waves of the Zakharov-Kuznetsov equation with critical speed, arXiv : 2004.10088v1 20 Apr 2020.

\bibitem{YRZ2019}
K. Yang, S. Roudenko and Y. Zhao, 
\emph{Blow-up dynamics in the mass super-critical NLS equations}, 
Phys. D, 396:47–69, 2019.


\bibitem{ZK} 
V.E. Zakharov and E.A. Kuznetsov, 
\emph{On three dimensional solitons}, 
Zhurnal Eksp. Teoret. Fiz, {\bf 66} (1974), 594--597 [in russian]; 
Sov. Phys JETP, vol. 39, no. 2 (1974), pp. 285--286.


\bibitem{NS2020}
{\tt http://www.mathphys.fr }    


\end{thebibliography}
\end{document}